%% file: chow_latex.tex
\documentclass[a4paper,11pt]{amsart}
\usepackage{amsthm,epsf,amssymb,amscd}

\def\oM{\overline{M}}
\def\uP{\underline{P}}

\def\bP{\Bbb P}

\def\dlog{operatorname{dlog}}

\def\Aff{\operatorname{Aff}}
\def\Sym{\operatorname{Sym}}
\def\Supp{\operatorname{Supp}}
\def\Hilb{\operatorname{Hilb}}
\def\Gr{\operatorname{Gr}}
\def\Chow{\operatorname{Chow}}

\def\tX{\tilde{X}}

\def\cO{\Cal O}
\def\tU{\tilde{U}}
\def\tB{\tilde{B}}

\def\Det{\operatorname{Det}}
\def\codim{\operatorname{codim}}
\def\rank{\operatorname{rank}}
\def\mult{\operatorname{mult}}
\def\bZ{\Bbb Z}
\def\bC{\Bbb C}
\def\bQ{\Bbb Q}
\def\bG{\Bbb G}
\def\bR{\Bbb R}
\def\bA{\Bbb A}
\def\bB{\Bbb B}
\def\bX{\Bbb X}
\def\oX{\overline{X}}
\def\obP{\bP^{\circ}}

\def\cB{\Cal B}
\def\cA{\Cal A}

\def\cS{\Cal S}
\def\cE{\Cal E}
\def\cF{\Cal F}
\def\tS{\tilde S}
\def\cI{\Cal I}
\def\cA{\Cal A}
\def\tcA{\tilde{\cA}}

\def\tX{\tilde{X}}
\def\tbP{\tilde{\bP}}

\def\ucI{\underline{\cI}}
\def\uP{{\underline{P}}}

\def\uI{\underline{I}}

\def\cT{\Cal T}
\def\cP{\Cal P}
\def\cB{\Cal B}
\def\tS{\tilde{S}}
\def\qp{\Cal P}
\def\cW{\Cal W}
\def\tqp{\tilde{\qp}}
\def\tB{\tilde{B}}
\def\tcB{\tilde{\cB}}
\def\oN{\overline{N}}

\def\tcA{\tilde{\cA}}
\def\cB{\Cal B}

\def\cX{\Cal X}
\def\cH{\Cal H}
\def\cM{\Cal M}

\def\Spec{\operatorname{Spec}}
\def\Proj{\operatorname{Proj}}
\def\PGL{\operatorname{PGL}}
\def\tM{\tilde{M}}

\def\oo{\overline{\Omega}}

\def\Stab{\operatorname{Stab}}
\def\aff{\Bbb A}
\def\cE{\Cal E}
\def\czE{\mathop{\cE}\limits^\circ}
\def\bS{\Bbb S}
\def\bB{\Bbb B}
\def\dlog{\operatorname{dlog}}

\def\qP{\Cal P}
\def\obS{\overline{\bS}}

\def\Res{\operatorname{Res}}
\def\Ord{\operatorname{Ord}}
\def\ord{\operatorname{ord}}

\def\ocB{\overline{\cB}}
\def\Star{\operatorname{Star}}
\def\oLambda{\overline{\Lambda}}

\def\osigma{\overline{\sigma}}
\def\bl{\operatorname{BL}}
\def\tqP{\tqp}
\def\ogamma{\overline{\gamma}}

\def\obS{\overline{\bS}}

\def\res{\Res}
\def\oZ{\overline{Z}}
\def\tPsi{\tilde{\Psi}}
\def\oR{\overline{R}}
\def\oK{\overline{K}}
\def\Sing{\operatorname{Sing}}
\def\Cal{\mathcal}

\newtheoremstyle{mystyle}{}{}{\itshape}{}{\scshape}{.}{ }{}
\theoremstyle{mystyle}

\swapnumbers
\newtheorem{Theorem}{Theorem}[section]

\newtheorem{ETheorem}[Theorem]{}
\newtheorem{Proposition}[Theorem]{Proposition}
\newtheorem{Lemma}[Theorem]{Lemma}
\newtheorem{Corollary}[Theorem]{Corollary}
\newtheorem{Claim}[Theorem]{Claim}
\newtheorem{Conjecture}[Theorem]{Conjecture}
\newtheorem{Question}[Theorem]{Question}

\newtheoremstyle{myreview}{}{}{}{}{\scshape}{.}{ }{}
\theoremstyle{myreview}

\newtheorem{Definition}[Theorem]{Definition}
\newtheorem{Example}[Theorem]{Example}

\newtheorem{Remark}[Theorem]{Remark}

\newtheorem{Notation}[Theorem]{Notation}
\newtheorem{Review}[Theorem]{}

\newcounter{et}[Theorem]
\def\cooltag{\tag{\arabic{section}.\arabic{Theorem}.\arabic{et}}\addtocounter{et}{1}}

\begin{document}

\title{Chow Quotients of Grassmannians II}
\author{ Sean Keel and Jenia Tevelev}

\address{Department of Mathematics, University of Texas at Austin,
Austin, Texas, 78712}
\email{keel@math.utexas.edu and tevelev@math.utexas.edu}

\maketitle

\section{Introduction and Statement of Results}\label{sean1}

\begin{Review}
Let
$\bP(r,n)$
be the space of ordered $n$-tuples of linear hyperplanes in~$\bP^{r-1}$. Let
$\obP(r,n) \subset \bP(r,n)$
be the open subset which are in linear general position.
$\PGL_r$ acts freely on
$\obP(r,n)$, let
$X(r,n)$ be the quotient.
$X(2,n)$ is usually denoted $M_{0,n}$. It has a compactification
$M_{0,n} \subset \oM_{0,n}$, due to Grothendieck and
Knudsen, with many remarkable properties:
\end{Review}

\begin{ETheorem}\label{sean1.0}{\sc Properties of $M_{0,n} \subset \oM_{0,n}$.}
\begin{enumerate}
\renewcommand\theenumi{\arabic{section}.\arabic{Theorem}.\arabic{enumi}}
\item\label{sean1.0.1} $\oM_{0,n}$ has a natural moduli interpretation,
namely it is the moduli space of stable $n$-pointed rational curves.
\item\label{sean1.0.2} Given power series $f_1(z),\dots,f_n(z)$ which
we think of as a one parameter family in $M_{0,n}$ one can
ask: What is the limiting stable $n$-pointed rational curve in
$\oM_{0,n}$ as $z\to0$ {\rm?} There is a beautiful answer,
due to Kapranov
\cite{Kapranov93a}, in terms of the Tits tree for
$\PGL_2$.
\item\label{sean1.0.3} $M_{0,n} \subset \oM_{0,n}$ has a natural Mori theoretic
meaning, namely it is the log canonical model, \cite{KM}.
In particular the pair $(\oM_{0,n},\partial \oM_{0,n})$
has log canonical singularities (a natural generalisation
of toroidal).
\end{enumerate}
\end{ETheorem}

In fact in \eqref{sean1.0.3} the pair has normal crossing, but we write the
weaker form as this is what there is a chance to generalize.
It is natural to wonder

\begin{Question}\label{sean1.1}
Is there a compactification
$X(r,n) \subset \oX(r,n)$ which satisfied any or all of
the properties of \eqref{sean1.0}\/{\rm?}
\end{Question}

\begin{Review}\textsc{Chow quotients of Grassmannians.}
There is an identification
$$
X(r,n) = G^0(r,n)/H,
$$
the so called Gelfand--Macpherson correspondence~\eqref{KHGKHGJ},
where $G(r,n)$ is the Grassmannian of $r$-planes in $\aff^n$,
$H=\bG_m^n$ is the standard diagonal torus and $G^0(r,n) \subset G(r,n)$ is
the open subset of $r$-planes which project isomorphically
onto $\aff^I$ for any subset $I \subset N$ with cardinality
$|I| =r$, where $N=\{1,\ldots,n\}$.
In \cite{Kapranov93},
Kapranov has introduced a natural compactification, the
so called Chow quotient
$$
X(r,n) = G^0(r,n)/H \subset G(r,n)//H := \oX(r,n).
$$
See \eqref{sean1.11} for a review of Chow quotients.
Here we note only that Kapranov defined a
natural flat family
\begin{equation}\label{sean1.2}
p:\, (\cS,\cB) \rightarrow \oX(r,n)\cooltag
\end{equation}
of pairs of schemes with boundary,
his so-called family of visible contours, generalising
the universal family over $\oM_{0,n}$, and Lafforgue in \cite{Lafforgue03}
gave a precise description of the fibres $(S,B)$, showing
in particular each pair has toroidal singularities
(throughout the paper $B$ will always indicate
a boundary, i.e. a Weil divisor, on a space clear from
context).
\end{Review}

$\oX(r,n)$ satisfies
the first two properties of \eqref{sean1.0}, but the third fails except
possibly in the 
exceptional cases $(2,n)$, $(3,6)$, $(3,7)$, $(3,8)$ (and those
obtained from these by a canonical duality). We conjecture
that in these cases the Chow quotient is indeed the log canonical
model, and speculate a relationship to the exceptional
root systems, see~\eqref{sean1.6}.
Let $\tX(r,n)\to\oX(r,n)$ be the normalization.
In~\eqref{KLJFKHGDKHD},
we follow Lafforgue and introduce the third modification $\oX_L(r,n)$
that sits between $\tX(r,n)$ and $\oX(r,n)$ (we distinguish 
between these spaces in the interests of precision, in fact the
minor differences will be for our purposes unimportant).

\begin{Review}\textsc{Moduli interpretation as in \eqref{sean1.0.1}.}
$\oX(r,n)$ is a natural moduli space of semi log
canonical pairs (the natural higher dimensional Mori theoretic
generalisation of stable pointed curves).
This is a recent result of Hacking, \cite{Hacking03}. We
observed the same result independently. Our proof, which
is based on \cite{Lafforgue03}, will appear elsewhere.
\end{Review}

\begin{Review}\textsc{Realization via Tits buildings as in \eqref{sean1.0.2}.}
Let $R = k[[z]]$ and $K$ its quotient field. Throughout
the paper $k$ is a fixed algebraically closed field. $V = k^r$
and $V_T = V\otimes_k T$ for
a $k$-algebra $T$. We write $\qP$ for projective
spaces of quotients (or equivalently, hyperplanes), 
$\bP$ for spaces of lines.
Begin with a collection
\begin{equation}\label{sean1.3}
\cF :=\{f_1,\dots,f_n \} \subset V_K \cooltag
\end{equation}
any $r$ of which are linearly independent, and thus
give a $K$-point of $X(r,n)$. We think of $\cF$ as
the equations of $n$ $1$-parameter families in $\bP^{r-1}$.
Following \cite{Kapranov93a} (where the case $r=2$ is treated) we
ask:

\begin{Question}
What is the limit as $z\to0$,
i.e. in the pullback of the family \eqref{sean1.2} along the associated
$R$-point of $\oX(r,n)$, what is the special fibre?
\end{Question}

We give a canonical solution, in terms of the Tits building
$\cB$. Here is a quick version, further details are given
at the end of this introduction, see \eqref{sean1.16}. Proofs and further
related results are given in \ref{sean4}--\ref{sean6}.
Recall $\cB$ is the
set of equivalence classes of
$R$-lattices in $V_K$ (i.e. free $R$-submodules of rank~$r$)
where $M$ and $N$ are equivalent if there exists $c \in K^*$
such that $M = cN$. A subset $Y \subset \cB$ is called {\it convex}
if it is closed under $R$-module sums, i.e. $[M_1],[M_2] \in Y$
implies $[M_1 + M_2] \in Y$.
For $Y \subset \cB$
we write $[Y]$ for its convex hull, which is finite if $Y$ is,
see \eqref{sean4.3.2}.

For a lattice $M$, and a non-trivial subset $\Theta\subset V_K$
(e.g.~an $R$-submodule or an element)
we can find unique $a > 0$ so that $z^a \Theta \subset M$,
$z^a \Theta\not \subset zM$. We define $\Theta^M := z^a \Theta \subset M$.
Let
$\Theta^{\oM} \subset \oM := M/zM$ denote the image of
the composition
$$
\Theta^M \subset M \twoheadrightarrow M/zM = \oM.
$$
We call $[M]$ {\it stable}
if $\cF^{\oM}$ contains $r+1$ elements in linear
general position. Let $\Stab \subset [\cF]$ be the
set of stable classes. $\Stab$ is finite, see \eqref{finitestable}.
\end{Review}

\begin{Definition}\label{LKUJHLKGKHG}
For a finite $Y \subset \cB$, let
$\bS_Y$ be the {\it join} of projective bundles $\qP(M)$,
$[M] \in Y$ -- i.e. fix one and take the closure of the graph
of the product of the birational maps from this projective
bundle to all the others.
\end{Definition}

Let $Y \subset \cB$ be any finite convex
subset containing $\Stab$ {\rm(}for example take the convex
hull $Y = [\Stab]${\rm)}. Let $\bB_i \subset \bS_Y$
be the closure of the hyperplane
$$
\{f_i = 0\} \subset \qP(V_K) \subset \bS_Y
$$
on the generic fibre of $p:\bS_Y \rightarrow \Spec(R)$. Let
$\bB = \sum \bB_i$, and $S_Y \subset \bS_Y$ the special fibre.
Let $V_n$ be the standard $k$-representation of $S_n$ (i.e.~elements
in $k^n$ whose coordinates sum to zero).

\begin{Theorem}\label{sean1.4}
$\bS_Y$, $\bB_i$ are non-singular and the divisor $S_Y + \bB$
has normal crossings. The $1$-forms $\dlog(f_i/f_j)$
define globally generating sections of the vector bundle 
$\Omega^1_{\bS/R}(\log \bB)$.
The image
of the associated map
$$
\bS_Y \rightarrow \Spec(R) \times G(r-1,n-1)
$$
is $\cS \rightarrow \Spec(R)$, the pullback of the family \eqref{sean1.2}
along the $R$-point of $\oX(r,n)$ defined by $\cF$. In particular
the relative log canonical bundle
$K_{\bS_Y} + \bB$ is relatively globally generated and big,
and $\bS_Y \rightarrow \cS$ is the relative minimal model, and
crepant.
\end{Theorem}

\begin{Review}\textsc{Historical Remark:}
By \cite[2.2]{Mustafin78}, if $Y$
is convex then $p:\bS_Y \rightarrow \Spec(R)$ is semi-stable, i.e.
$\bS_Y$ is non-singular, and the closed fibre $S_Y$ has smooth
irreducible components and normal crossings.
Mustafin remarks after the
proof that the join probably represents some
natural functor but he prefers the explicit join construction.
This functor (see \eqref{sean1.14}) is introduced in Definition~4 of \cite{Faltings01},
and Faltings proves it represents the join.
Faltings attributes the functorial description to Deligne,
to whom Mustafin also refers. We will refer to $\bS_Y$
as the Deligne scheme, or Deligne functor.
\end{Review}

For the definition and basic properties of bundles of relative log
differentials see \ref{sean9}.
We note that the crepant semi-stable model
$(\bS_Y,\bB)$ is in many ways preferable to its minimal model
$\cS \rightarrow \Spec(R)$. For example dropping the last
hyperplane induces a natural regular birational map
$$
\bS_{[\Stab(\cF)]} \rightarrow \bS_{[\Stab(\cF')]}
$$
for $\cF' = \cF \setminus \{f_n\}$, but for $r \geq 3$ the associated
rational map between minimal models is not in general regular.
There are examples where regularity fails already with $(r,n) = (3,5)$.

The special fibres $(S_Y,B)$ and $(S,B)$ can be canonically
described in terms of the membrane $[\cF] \subset \cB$. We turn
to this in \eqref{sean1.16} below, but wish first to discuss \eqref{sean1.0.3}:

\begin{Review}\textsc{Log canonical model as in  \eqref{sean1.0.3}.}\label{sean1.7}
Let $M$ be a smooth variety over the complex numbers,
and let $M \subset \tM$ be a compactification, with
normal crossing divisorial boundary, $B$. The vector spaces
$$
H^0(\tM,\omega_{\tM}(B)^{\otimes m})
$$
turn out to depend only on $M$, and so give a
canonical rational map, the so called $m$-plurilogcanonical
map, to projective space. The finite generation conjecture
of Mori theory implies that
if for some $m$ the map is an immersion,
then for sufficiently large $m$, the closure of the
image gives a compactification $M \subset \oM$
independent of $m$, and with boundary, $B$ having
nice singularities, namely $K_{\oM} + B$ is log
canonical. We do not recall here the definition
of log canonical (see e.g. \cite{KMM87}) but note
that one can think of it as a weakening of toroidal,
the pair of a toric variety with its boundary being one
example.

The initial motivation for this paper was
the elementary observation, \eqref{sean1.8}, that
$X(r,n)$ is minimal of log general type (its first log
canonical map is a regular immersion) and thus there is
(conjecturally) a natural Mori theoretic compactification,
the log canonical model. It is natural to wonder:

\begin{Question}\label{sean1.9}
What is the log canonical
model $X(r,n) \subset \oX_{lc}(r,n)$?
\end{Question}

We believe this compactification is
of compelling interest, as it gives a birational
model with reasonable boundary singularities of a compactification
of $X(r,n)$ whose boundary components meet in absolutely
arbitrary ways. In particular, were $\oX(3,n)$ the log canonical
model, it would give something like a canonical resolution of
all singularities, see \eqref{sean1.10}.
Unfortunately, the two compactifications do not in general agree:
\end{Review}

\begin{Theorem}\label{sean1.5}\label{sean3.15}\label{sean3.18}
$\tX(3,n)$ with its boundary fails to be
log canonical for $n \geq 9$
{\rm(}for $n \geq 7$ in characteristic $2${\rm)}.
$\tX(4,n)$ is not log canonical for $n\ge8$.
\end{Theorem}

We prove \eqref{sean1.5} in \ref{sean3}.
Moreover, \eqref{sean3.17} shows that in general
the pair $(\oX(3,n),B)$ has arbitrary singularities. We note
\eqref{sean1.5} is at variance with the hope expressed in
\cite{Hacking03}, and one that we ourselves for a long time
harbored, that the pair has toroidal singularities. Faltings
and Lafforgue, \cite{Faltings01}, \cite{Lafforgue99} expressed
the same hope for their compactification of
$\PGL_r^n/\PGL_r$ (which is itself a Chow quotient in a natural
way), but
Lafforgue has shown this hope was false as well,
\cite[3.28]{Lafforgue03}.

In the positive direction we speculate the two
agree in the cases that remain:

\begin{Conjecture}\label{sean1.6}
$X(r,n) \subset \oX(r,n)$ is
the log canonical model precisely in the cases
$$
(2,n),\ (3,6),\ (3,7),\ (3,8)
$$
and those obtained from these by the
canonical duality \cite[2.3]{Kapranov93}
$$\oX(r,n) = \oX(n-r,n).$$
Moreover in these cases the pair $(\oX(r,n),B)$ has toroidal
singularities.
\end{Conjecture}

The numbers in \eqref{sean1.6} are of course very suggestive and
it is natural to wonder if there is a connection
with the
exceptional root systems
$$
D_n,\ E_6,\ E_7,\ E_8.
$$

\begin{Review}\textsc{Equations and syzygies of $\oM_{0,n}$.}
The closure of $X(2,n)$ in the first log canonical
immersion is $\oX(2,n)=\oM_{0,n}$
(in particular, \eqref{sean1.6} holds for $\oX(2,n)$)
and the equations for the embedding are nice:

\begin{Theorem}\label{sean2.17}
Let $\kappa := \omega_{\oM_{0,N}}(B)$.
Let
$$
N_3 \subset N_4 \subset \dots \subset N_{n} = N
$$
be a flag of subsets of $N$, with $|N_j| = j$.
There is a canonical identification
$$
H^0(\oM_{0,N},\omega(B)) = \bigotimes_{j=3}^{n-1} V_{N_{j}}
$$
and over the complex numbers a canonical identification
$$
H^0(\oM_{0,N},\kappa) = H^{n-3}(M_{0,N},\bC).
$$
$\kappa$ is very ample, the embedding factors through
the Segre embedding
$$
\oM_{0,N} \subset \bP(V_{N_3}^\vee) \times
\bP(V_{N_4}^\vee) \dots \times \bP(V_{N_{n-1}}^\vee)
\subset \bP(H^0(\oM_{0,N},\kappa)^\vee)
$$
and $\oM_{0,N}$ is the scheme theoretic intersection of the
Segre embeddings over flags of subsets, i.e.
$$
\oM_{0,N} = \bigcap_{{N_3 \subset N_4 \dots \subset N}}
\bP(V_{N_3}^\vee) \times \bP(V_{N_4}^\vee) \times \bP(V_{N_{n-1}}^\vee)
\subset \bP(H^0(\oM_{0,N},\kappa)^\vee).
$$

$$
\Sym(H^0(\oM_{0,N},\kappa)) \rightarrow \bigoplus_{n \geq 0}
H^0(\oM_{0,N},\kappa^{\otimes n})
$$
is surjective, and the kernel is generated by quadrics, and
the syzygies among the quadrics are generated by linear
syzygies.
\end{Theorem}

Above $V_{N_i}$ is the standard $k$-representation of
the symmetric group $S_{N_i}$.
The proof of \eqref{sean2.17} will appear elsewhere.
\eqref{sean2.17} implies in particular that the compactification
$M_{0,n} \subset \oM_{0,n}$ can be recovered in a canonical way
from the system of $S_i$ modules $H^{i-3}(M_{0,i},\bZ)$ together
with the pullback maps between them (for the fibrations given
by dropping points), i.e.~can be recovered
canonically from the Lie operad, see \cite{Getzler}.
\end{Review}

\begin{Review}\textsc{The Chow quotient of $G(3,6)$.}
We have proven \eqref{sean1.6}
also for $\oX(3,6)$. This space is very interesting. For example:
There is a natural map $\oX(r,n+1) \rightarrow \oX(r,n)$ given
by dropping the last hyperplane. For $\oX(2,n)$ it is
well known that the map is flat, and
canonically identified with the universal family \eqref{sean1.2}.
For $\oX(3,5)$ it is again flat (this fails for $\oX(3,n)$,
$n \geq 6$), and a natural universal family. However it is
not the family \eqref{sean1.2},
but rather Lafforgue's
analogous family (he defines such a family beginning with
any configuration of hyperplanes, see \ref{sean2}) for the
configuration of $10$ lines which is dual to the configuration
of $5$ general lines.
There are $15$ irreducible components of the boundary
$B \subset \oX(3,6)$
which surject onto $\oX(3,5)$. If we let $\Gamma$ be
their union, then
$$
(\oX(3,6),\Gamma) \rightarrow \oX(3,5)
$$
gives a flat family of pairs, compactifying the family
of pairs $(S,B)$ for $S$
a del Pezzo surface of degree $4$, and $B$ a union of some of its
$-1$ curves. A detailed study of $\oX(3,6)$, including proofs
of all the claims of this paragraph, will appear elsewhere.
\end{Review}

\begin{Review}
Using results of \cite{Lafforgue03} we give a cohomological criterion
under which $\oX(r,n)$ will be a log minimal model, see \eqref{sean2.16}.
This we expect will apply in the cases of \eqref{sean1.6}. Note the
statement
$M \subset \oM$ is the log canonical model has two, in general
entirely independent, parts: First a singularity statement,
$(\oM,B)$ is log canonical (morally,
toroidal), and second a positivity statement,
$K_{\oM} + B$ is ample. However by \eqref{sean1.2} and \eqref{sean2.16}
in the case of
of $X(r,n) \subset \oX(r,n)$ it turns out that whenever the first happens,
the second comes for free. See \eqref{sean2.16} for the precise statement.
\end{Review}

\begin{Review}\label{sean1.10}\textsc{Mnev's universality theorem.}
The boundary
$$\bP(r,n) \setminus \obP(r,n)$$
is a union of $\binom{n}{r}$ Weil divisors. The
components have only mild singularities, however they
meet in very complicated ways:
Let $Y$ be an affine scheme of
finite type over $\Spec(\bZ)$. By \cite[1.8]{Lafforgue03},
there are integers $n,m$
and an open subset
$$
U \subset Y \times \aff^m
$$
such that the projection $U \rightarrow Y$ is surjective,
and $U$ is isomorphic to a boundary stratum of
$\bP(3,n)$ (a boundary stratum for a divisorial boundary means the
locally closed subset of points which lie in each of a prescribed
subset of the irreducible components, but no others).

An analogous statement holds for the boundary
of $X(3,n)$, say in any of its smooth GIT quotient
compactifications.
Our \eqref{sean3.17} suggests the singularities
of the pair $(\oX(3,n),B)$ are also in general arbitrary scheme-theoretically.
Now by \eqref{sean1.8} and the finite generation
conjecture of Mori theory,
$X(3,n) \subset \oX_{lc}(3,n)$ will give a
canonical compactification in which the boundary has
mild (i.e. log canonical) singularities. $\oX(r,n)$ maps
to all the GIT quotients. We do not know whether or not this
will hold for $\oX_{lc}(3,n)$ -- if it does $\oX_{lc}$ would
give an absolutely canonical way of (partially) resolving
a boundary whose strata include all possible singularities.
\end{Review}

Now we return to the Tits building,
to give a canonical description of special fibres.

\begin{Definition}
We define the {\it membrane}, $[\cF] \subset \cB$
to be classes of lattices which have a basis
given by scalar multiples of some $r$ elements from $\cF$,
or equivalently, such that the limits $\cF^{\oM}$ contain
a basis of $\oM$.
\end{Definition}

\begin{Review}\label{sean1.16}\textsc{Special fibers.}
The building $\cB$ is a simplicial complex of dimension $r-1$:
We say $[M],[N] \in \cB$
are {\it incident} if we can choose representatives so
$$
zM \subset N \subset M
$$
(the relation is easily seen to be symmetric). Points
$[M_1],\dots,[M_m]$ span an $(m-1)$-simplex iff they are
pairwise incident, which holds iff we can
choose representatives so
$$
zM_m = M_0 \subset M_1 \subset \dots \subset M_m.
$$
By scaling we can put any of the $M_i$ in the position
of $M_m$ but the cyclic ordering among them is intrinsic.

Now take a convex subset
$$Y \subset [\cF] \subset \cB$$
(not necessarily finite).
Canonically associated to each $(m-1)$-simplex
$\sigma \subset Y$ as above, is a smooth projective variety
$\tqP(\osigma)$,
$$
\tqP(\osigma) = \prod_{m \geq i \geq 1} \tqP(M_{i}/M_{i-1})
$$
where $\tqP(M_{i}/M_{i-1})$ is a certain iterated blowup
of the projective space of quotients $\qP(M_i/M_{i-1})$
along smooth centers. For precise details see \eqref{sean5.6}, \eqref{sean5.7}.
There are canonical compatible closed embedding
$\tqP(\osigma) \subset \tqP(\ogamma)$ for simplicies
$\gamma \subset \sigma \subset Y$. Finally there
is a canonical scheme $S_Y$, with irreducible components
$\tqP(\oM), [M] \in Y$, such that for a subset
$\sigma \subset Y$, the $\tqP(\oM)$, $[M] \in \sigma$
have common intersection iff $\sigma$ is a simplex,
and in that case the intersection is
$\tqP(\sigma)$, e.g. $\tqP(\oM)$ and $\tqP(\oN)$
for lattices $[M],[N] \in Y$ meet iff they span a
$1$-simplex, $\sigma$, and in that case they are glued along
the common smooth divisor $\tqP(\osigma)$. In particular,
$S_Y$ has smooth components and normal crossings.
When $Y$ is finite, $S_Y$ is the special fiber of
$\bS_Y$. It carries Cartier boundary
divisors $B_i \subset S_Y$ for each $i \in N$. These are described
as follows: $B_i$ has smooth irreducible components, and
$\sum B_i$ has normal crossings. $B_i$
has a component on
$\tqP(\oM) \subset S_Y$, $[M] \in Y$ iff
the lattice $[M + z^{-1} f_i^{M}] \in \cB$ is not in $Y$.
In this case the component is the strict transform
of the hypersurface $\{f_i^{\oM} = 0\} \subset \qP(\oM)$.
The limit
variety $(S,B)$ (i.e. the fibre
of $(\cS,\cB)$ over the image of the closed point
of $R$) is the $K_{S_Y} + B$ minimal model, see~\eqref{sean1.18}.
\end{Review}

\begin{Review}\textsc{Bubble space.}
For $Y = [\cF]$, write $S_{\infty} = S_Y$. For
$r=2$, $S_{\infty}$ is the scheme constructed
in \cite{Kapranov93a} -- it is a tree of rational curves
with countably many components such that each component intersects
at least two others. $S_{\infty}$ has no boundary, its canonical
linear series $|K_{S_{\infty}}|$ is
globally generated and the image of the associated map is
again $S$. More precisely:
\end{Review}

\begin{Theorem}\label{sean1.17}
$S_{\infty}$ has smooth projective components
and normal crossings. It carries a natural vector bundle
$\Omega^1(\log)$, with determinant $\omega_{S_\infty}$.  For each finite
convex subset
$$
\Stab \subset Y \subset [\cF]
$$
there is a canonical
regular surjection
$$
p:\, S_{\infty} \rightarrow S_Y
$$
and a canonical isomorphism
$$
p^*(\Omega^1(\log B)) \rightarrow \Omega^1(\log).
$$
Given a closed point $x \in S_{\infty}$ there exists
a $Y$ so that $p$ is an isomorphism in a neighborhood
of $x$.

The differential forms $\dlog(f/g)$, $f,g \in \cF$,
induce a canonical inclusion
$$
V_n \subset H^0(S_{\infty},\Omega^1(\log)).
$$
These sections generate the bundle. The associated map
$$S_{\infty} \rightarrow G(r-1,n-1)$$ factors through
$S_Y$ and the image is the limit variety $S$.
\end{Theorem}

The picture illustrates the case $r=2$:

\medskip
\epsfxsize=0.9\textwidth
\centerline{\epsfbox{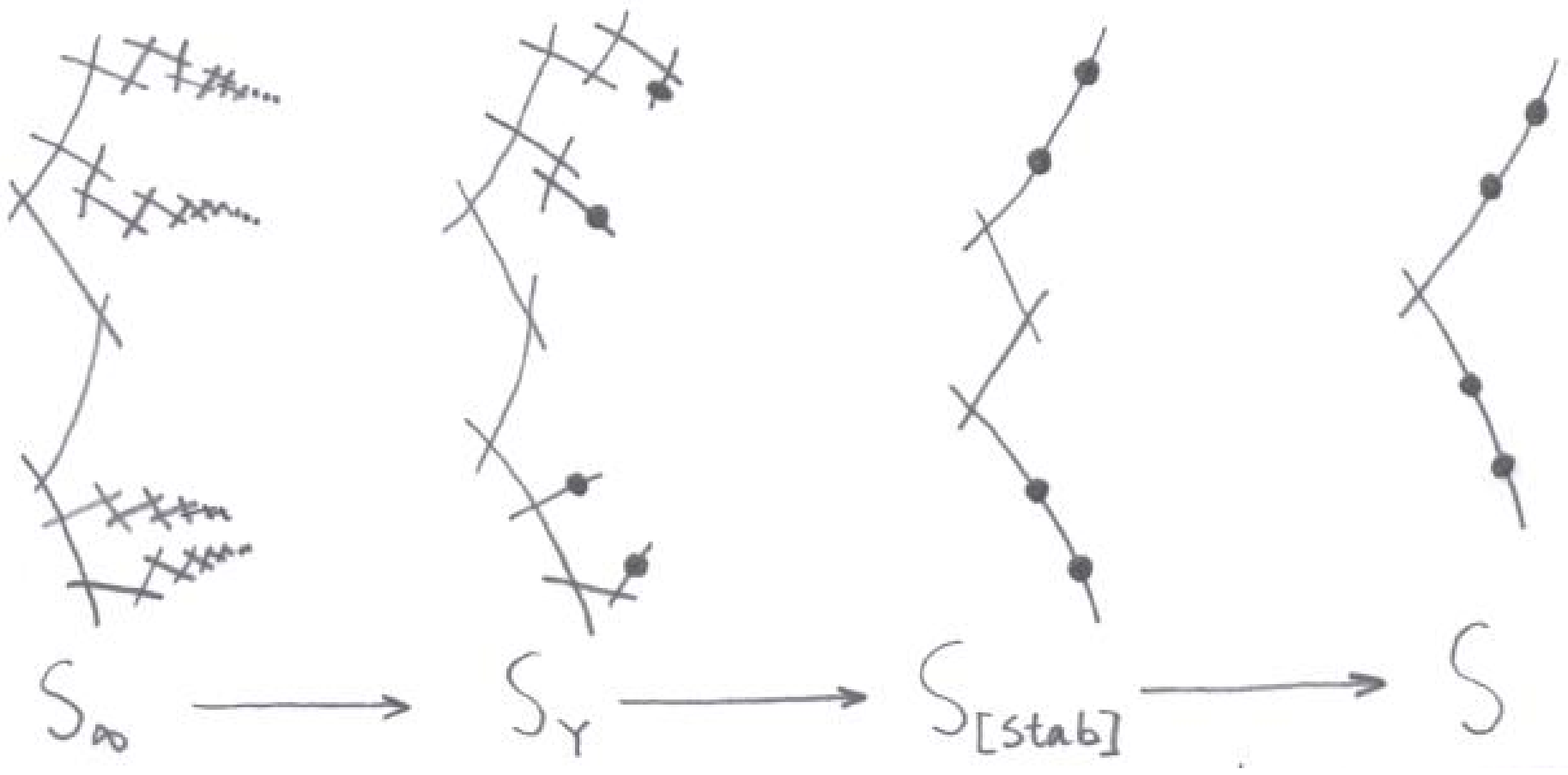}}

\begin{Review}\label{sean1.18}\textsc{The limit variety}
$(S,B)$ can also be canonically recovered
from the membrane: For each $[M] \in [\cF]$,
$\tqP(\oM)$ has a canonical normal crossing boundary, $B$,
the union of the divisors $\tqP(\osigma) \subset \tqP(\oM)$ over
$1$-simplicies $[M] \in \sigma \subset [\cF]$. The
rational differential forms $\dlog(f^{\oM}/g^{\oM})$ on
$\tqP(\oM)$ have log poles along $B$, and so define
canonical sections of $\Omega^1(\log B)$. These sections
generate the bundle. In particular their wedges generate
$K_{\tqP(\oM)}+B$.

\begin{Definition}\label{sean1.15}
We call a configuration of hyperplanes {\it GIT stable}
if its group of automorphisms is trivial.
We call $[M] \in [\cF]$ GIT stable if the
configuration of limiting hyperplanes $\cF^{\oM}$ is GIT stable.
Of course stable implies GIT stable. For $r \leq 3$ they
are the same, see \ref{sean7.1}, but they are in general different. There
are only finitely many GIT stable equivalence classes, see \ref{sean6.16}.
$[M]$ is GIT stable if and only if $K_{\tqP(\oM)} + B$ is big.
\end{Definition}

The irreducible components of
$S$ are the $(K_{\tqP(\oM)}+B)$-minimal models of
$\tqP(\oM)$ for GIT stable  $\oM$.
The minimal model can be constructed
as follows: Let $U_{\oM} \subset \qP(\oM)$ be
the complement to the union of hyperplanes. Equivalently,
$U_{\oM} = \tqP(\oM) \setminus B$. The (regular) differential forms
$\dlog(f^{\oM}/g^{\oM})$ generate the cotangent bundle of
$U_{\oM}$, the associated map to $G(r-1,n-1)$ is an immersion,
and the corresponding irreducible component of $S$ is the
closure of $U_{\oM} \subset G(r-1,n-1)$.
\end{Review}

\begin{Review}\label{sean1.19}\textsc{Illustration.}
Let us look at the first non-trivial example for $r=3$.
$$
\cF = \{f_1, f_2, f_3,\quad f_4 = f_1+f_2+f_3,\quad  f_5 = z^{-1}f_1+f_2+f_3\}
$$
for $f_1,f_2,f_3$ the standard (constant) basis of $k^3 \subset K^3$.
In this case there are two stable lattices,
$$
M_1 = R f_1 + Rf_2 + R f_3\quad\hbox{\rm and}\quad
M_2 = R z^{-1} f_1 + R f_2 + R f_3.
$$

\epsfxsize=0.8\textwidth
\centerline{\epsfbox{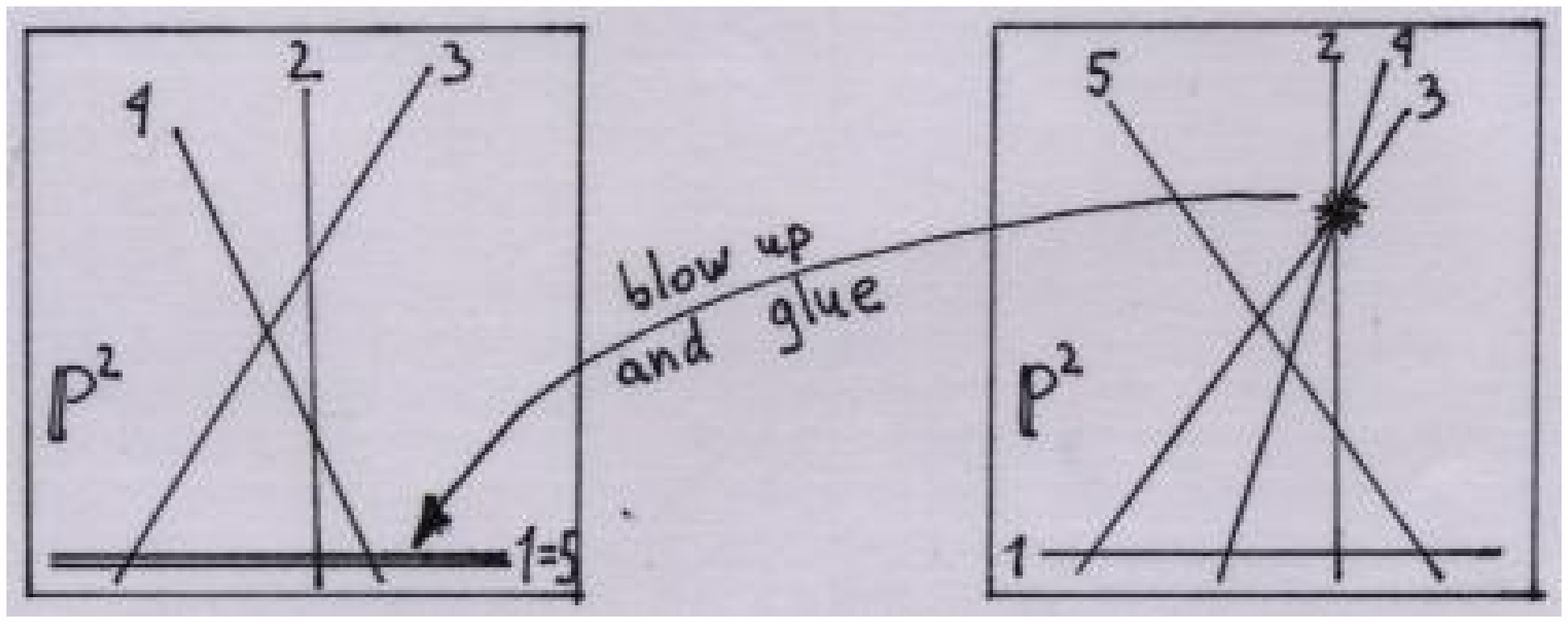}}
\medskip
The picture
illustrates the limit surface $(S,B)$.
Note
$zM_2 \subset M_1 \subset M_2$, so the
stable lattices in this case form a $1$-simplex, $\sigma$,
in
particular are already convex. So we can take
$Y = \sigma = \Stab$.
The two pictures are the configurations of limit lines.
The components of $S$ are
$\tqP(\oM_1) = \qP(\oM_1)$, and $\tqP(\oM_2)$, the blowup
of $\qP(\oM_2)$ at the intersection
point  $f_2^{\oM_3}= f_3^{\oM_3}=0$. The two components are
glued along $\tqP(\osigma) = \qP(M_2/M_1)$,
which embeds in $\tqP(\oM_1)$
as the line $f_1^{\oM_1} =0$, and in $\tqP(\oM_2)$ as
the exceptional curve.

Unfortunately we can't draw the membrane:
if $n\ge5$ then no membrane of $\oX(3,n)$ can be
embedded in $\bR^3$ without self-intersections.
\end{Review}

\begin{Review}\textsc{\textsc Relation to tropical algebraic geometry.}
There are several connections between this work and
tropical algebraic geometry:
$[\cF]$ is naturally
homeomorphic to the tropicalisation of the $r$-dimensional subspace in
$K^n$ defined by the rows of the matrix with columns
the $f_i$, see \eqref{sean4.5.1}. Further we observe, \eqref{sean4.5.2} a natural
generalisation of Kapranov's family \eqref{sean1.2} which might give information
about an arbitrary tropical variety. It may also
help explain an interesting correspondence, mysterious to us at present:
We have observed that the incident
combinatorics of the boundary of $\oX(r,n)$ are
encoded in the tropical Grassmannian, introduced
in \cite{SS}, at
least in the cases $\oX(2,n)$ and $\oX(3,6)$, which
are the only cases in which the tropical Grassmannian, or
the boundary strata, have been computed.
\end{Review}

We thank B. Hassett, J. McKernan,
and F. Ambro for help understanding material related
to the paper. W. Fulton pointed out to us
a serious error in an earlier version of the introduction.
M. Olsson helped us a great deal with log structures,
in particular we learned the
construction \eqref{sean2.11} from him.
We would like to especially thank L. Lafforgue for a series of
detailed email tutorials on \cite{Lafforgue03} and
M. Kapranov for
several illuminating discussions, and in particular for posing
to us the question of whether his \ref{sean1.0.2}
could be generalized to higher dimensions. The first author
was partially supported by NSF grant DMS-$9988874$.

\section{Various toric quotients of the Grassmannian}\label{sean2}

\begin{Review}\textsc{Chow variety \cite{Barlet}.}
Let $\Chow_{k,d}(\bP^{n})$ be the set of all $k$-dimensional algebraic cycles
of degree $d$ in $\bP^n$.
There is a canonical embedding
$$\Chow_{k,d}(\bP^{n})\subset\bP({\Cal V}), \quad X\mapsto R_X,$$
where $R_X$ is a Chow form of $X$. Here ${\Cal V}=H^0(G(n-k,n+1),\cO(d))^\vee$.
The image is Zariski closed, so
$\Chow_{k,d}(\bP^{n})$ is a projective variety with a canonical Chow polarization.

Now if $X\subset\bP^{n}$ is any projective subvariety and $\delta\in H_{2k}(X,\bZ)$
then the set $\Chow_\delta(X)$ of algebraic cycles in $X$ with the homology class $\delta$
is a Zariski closed subset of $\Chow_{k,d}(\bP^{n})$, where
$d\in H_{2k}(\bP^{n},\bZ)$ is the image of $\delta$.
The resulting structure of the algebraic variety on $\Chow_\delta(X)$ does not
depend on the projective embedding.
\end{Review}

\begin{Review}\label{sean1.11}\textsc{Chow quotients \cite{KSZ91}.}
Let $H$ be an algebraic group acting on a projective variety $G$.
Let $G^0 \subset G$ be a (sufficiently) generic open $H$-invariant subset.
In particular, all orbit closures $\overline{Hx}$, $x\in G^0$ have the same homology class~$\delta$.
There is a natural map
$$
G^0/H \rightarrow \Chow_\delta(G),\quad
x\mapsto\overline{Hx}.
$$
The Chow quotient $G//H$
is the closure of the image of this map.
There is also a parallel theory of Hilbert quotients $G///H$
when one takes the closure of $G^0/H$ in the Hilbert scheme of $G$.
\end{Review}

\begin{Review}\textsc{Chow quotients of projective spaces \cite{KSZ91,GKZ94}.}\label{HGPOU}
Let $H$ be an algebraic torus with the character lattice $\bX$.
Let $P\subset\bX_\bR$ be a convex polytope with vertices in~$\bX$.
We will denote vertices of $P$ by the same letter.
Let $V$ be a $k$-vector space with a basis $\{z_p\,|\,p\in P\}$.
The torus $H$ acts on~$V$ by formula $h\cdot z_p=p(h)z_p$.
We are going to describe $\bP//H$, where $\bP=\bP(V)$.
For any $S\subset\bP$, let $\Supp(S)\subset P$ be the set of coordinates
that don't vanish on $S$.
Let $\bP^0=\{p\in\bP\,|\,\Supp(p)=P\}$.

Take the big torus $\cH=\bG_m^P$ with its obvious action on~$V$
(so $P$ is identified with the set of ``coordinate'' characters of $\cH$).
We can assume without loss of generality that $H\subset\cH$.
This is equivalent to $\langle P\rangle_\bZ=\bX$, where
for any $S\subset\bX$ we denote by $\langle S\rangle_\bZ$ the minimal sublattice
containing $S$.

Actions of $H$ and $\cH$ on $\bP$ commute, therefore $\cH$ acts on $\bP//H$.
Moreover, since all points in $\bP^0$ are $\cH$-equivalent, $\cH$ has an open orbit $\bP^0/H\subset\bP//H$.
So $\bP//H$ is a projective $\cH$-toric variety
with the canonical $\cH$-equivariant Chow polarization.
By a toric variety we mean a variety with an action of a torus having a dense open orbit.
We don't assume that the action is effective or that the variety is normal.

Let $\Psi:\,\bP{//}_nH\to\bP//H$ be the normalization.

\begin{Remark}
$\Psi$ is bijective on the set of orbits
(this is true for any projective toric variety with an equivariant 
polarization)
but (as far as we know) is not always a bijective map.
We study this issue in more detail at the end of this section, 
see~\eqref{JKFKJGF}.
\end{Remark}

The fan $\cF(P)$  of $\bP{//}_nH$ can be described as follows.
A triangulation $T$ of~$P$ (with all vertices in the set of vertices of $P$) is called coherent if there exists a concave
piecewise affine function on $P$ whose domains of affinity are precisely maximal simplices of~$T$.
It gives rise to a polyhedral cone $C(T)\subset\bR^P$ of the maximal dimension.
Namely, $C(T)$ consists of all functions $\psi:\,P\to\bR$ such that $\psi_T:\,P\to\bR$
is concave, where $\psi_T$ is given by affinely interpolating $\psi$ inside each simplex of $T$.
Cones $C(T)$ (and their faces) for various $T$ give a complete fan $\cF(P)$.
Lower-dimensional faces of $\cF(P)$
correspond to (coherent) polyhedral decompositions $\uP$ of $P$.
More precisely, $C(\uP)$ is the set of concave functions
affine on each polytope of~$\uP$.

We will on ocassion abuse notation and refer to the collection
of maximal dimensional polytopes of a polyhedral decomposition
as a polyhedral composition itself.

We have the orbit decomposition
$$\bP//H=\bigsqcup_\uP\,(\bP//H)_\uP$$
(and a similar one for $\bP{//}_nH$)
indexed by polyhedral decompositions.
A cycle $X\in(\bP//H)_\uP$ is the union of toric orbits
with multiplicities, moreover
\begin{equation}\label{brokentoric}
X=\sum_{P'\in\uP}m_{P'}X_{P'},\quad
\Supp(X_{P'})=P',\quad
m_{P'}=[\bX:\langle P'\rangle_\bZ].
\cooltag
\end{equation}
If $m_{P'}=1$ for any $P'\in\uP$ then we say that $X$ is a {\em broken toric variety}.
If
$$\langle P''\rangle_\bZ=\bX\cap\langle P''\rangle_\bQ$$
for any face $P''$ of a polytope $P'\in\uP$
then we call $\uP$ unimodular.
\end{Review}

\begin{Review}\textsc{Hypersymplex.}
Let $H=\bG_m^n$ be the standard torus acting on $\wedge^r k^n$.
The weights $e_{i_1}+\ldots+e_{i_r}\in\bR^n$ are the vertices of the hypersymplex
$$
\Delta(r,n) := \left\{(x_1,\dots,x_n) \in \bR^n\ \big|\ \sum x_i = r,\ 1 \geq x_i \geq 0\right\}.
$$
$\Delta(r,n)$ has $2n$ faces $\{x_i=0\}$ and $\{x_i=1\}$.
The Pl\"ucker embedding $G(r,n)\subset\bP(\wedge^r k^n)$
induces a closed embedding $G(r,n)//H\subset\bP(\wedge^r k^n)//H$.
\end{Review}

\begin{Review}\textsc{Lafforgue's variety $\cA$.}\label{KJHGKHGJG}
If $x\in G(r,n)\subset\bP(\wedge^r k^n)$ then
a convex hull of $\Supp(x)\subset\Delta(r,n)$
is always a so-called matroid polytope (for the definition see~\eqref{KJHGKG}).
Lafforgue defines his varieties for an arbitrary fixed matroid polytope $P$.
Let
$$\bP^P=\{x\in\bP\,|\,\Supp(x)\subset P\},\quad
\bP^{P, 0}=\{x\in\bP^P\,|\,\Supp(x)=P\}$$
The locally closed subscheme
$$G^{P, 0}(r,n)=G(r,n)\cap\bP^{P, 0}$$
is called a thin Schubert cell.
Of course, $G^{\Delta(r,n), 0}(r,n)=G^0(r,n)$.
Lafforgue's scheme $\oo^P$ (see~\ref{jhfjhgfjfjfd})
is a compactification of $G^{P, 0}(r,n)/H$.
In almost all our applications, $P=\Delta(r,n)$
and so we adopt the following Notational Convention throughout
the paper: If we drop the polytope $P$ from
notation, it is assumed to be $\Delta(r,n)$, for a pair $(r,n)$ clear
from context.

In \cite[2.1]{Lafforgue03} Lafforgue defines a subfan of $\cF(P)$ whose cones
are in one to one correspondence with matroid decompositions $\uP$ of $P$ (i.e.~tilings of~$P$ by matroid polytopes).
This is a fan because a polyhedral decomposition coarser than a matroid decomposition
is a matroid decomposition (moreover,
if a convex polytope $Q\subset\Delta(r,n)$ admits a tiling
by matroid polytopes then $Q$ itself is a matroid polytope, see~\cite{Lafforgue99}).
The associated toric variety is denoted~$\cA^P$.
Just by definition, $\cA^P$ is the toric open
subset in the normalization of the Chow quotient:
\begin{equation}\label{sean2.1}
\cA^P\subset \bP^P{//}_nH. \cooltag
\end{equation}
Orbits in $\cA^P$ correspond to matroid decompositions.
Notice that the action of $\bG_m^P$ on $\cA^P$ is not effective:
the kernel $(\bG_m^P)_{\emptyset} \subset \bG_m^P$ is
the subtorus of affine maps $P\to\bG_m$.
Let  $\cA^P_{\emptyset} := \bG_m^P/(\bG_m^P)_{\emptyset}$.

For any face $Q$ of $P$, Lafforgue defines a natural {\em face map} of toric varieties
$\cA^P\to\cA^Q$. The corresponding map of fans is given by the restriction
of piecewise affine functions from $P$ to $Q$.
In particular, the image of the orbit $\cA^P_\uP$ is $\cA^Q_{\uP'}$,
where $\uP'$ is the matroid decomposition of $Q$ obtained by intersecting
polytopes in $\uP$ with  $Q$.
\end{Review}

\begin{Review}\textsc{Lafforgue's variety $\tcA$.}
Lafforgue introduces a second normal toric variety~$\tcA^P$ for the torus
$\tcA^P_{\emptyset} := \bG_m^P/\bG_m$ and a map of toric varieties
\begin{equation}\label{sean2.2}
\tcA^P \rightarrow \cA^P \cooltag
\end{equation}
extending the natural quotient map
$\tcA^P_{\emptyset} \twoheadrightarrow \cA^P_{\emptyset}$.

The torus orbits of $\tcA^P$ are in one to one correspondence
with $(\uP,P')$, for $\uP$ a matroid decomposition, and
$P' \in \uP$ one of the matroid polytopes.

By \cite[Proposition IV.3]{Lafforgue03}, \eqref{sean2.2} is projective and flat,
with geometrically reduced fibres, and there exists a natural equivariant
closed embedding
\begin{equation}\label{KJYFKJTF}
\tcA^P  \subset \cA^P \times \bP^P. \cooltag
\end{equation}

The fibre of \eqref{sean2.2} over a closed point of $\cA^P_\uP=(\bP^P{//}_nH)_\uP$
is a broken toric variety \eqref{brokentoric} in $\bP^P$.
All multiplicities are equal to~$1$ because of the following fundamental observation~\cite{GMS87}:
\begin{equation}\label{KFJHGDFJFDJKHGD}
\hbox{\em any matroid decomposition is unimodular.}\cooltag
\end{equation}
In fact, \eqref{sean2.2} is
the pullback of the universal Chow family over the Chow quotient
$\bP^P//H$
along the map $\cA^P\subset \bP^P{//}_nH\to \bP^P//H$.

For each maximal face $P'$ of $P$, the pair $(\emptyset,P')$,
where $\emptyset$ denotes the trivial decomposition (just
$P$ and its faces), corresponds to an
irreducible boundary divisor of $\tcA^P$. Denote
the union of these boundary divisors as
$\tcB^P \subset \tcA^P$. In the case
of $P = \Delta(r,n)$ there are $2n$ such boundary divisors,
corresponding to the maximal faces $\{x_i =0\}$, $\{x_i =1\}$ of
$\Delta(r,n)$. We indicate by $\tcB_i$ the divisor corresponding
to $\{x_i =1\}$.
Boundary divisors of $\tcA^P$ induce boundary divisors $B$ on fibres of
\eqref{sean2.2} for each maximal face of $P$. For $P = \Delta(r,n)$ we write $B_i$ for the
divisor corresponding to $\tcB_i$.
\end{Review}

\begin{Review}\label{jhfjhgfjfjfd}\textsc{Lafforgue's compactification $\oo$.}
Next we consider Lafforgue's main object, $\oo^{P}$, which
we consider only in the case $P = \Delta(r,n)$. We use
a different construction from his, as it is a quicker way of
describing the scheme structure -- $\oo$ is the subscheme
of $\cA$ over which the fibres of \eqref{sean2.2} are contained in
$G(r,n)$.                  
\end{Review}


\begin{Proposition}\label{sean2.4}
The Lafforgue space
$\oo \subset \cA$ is $\varphi^{-1}(\Hilb(G(r,n)))$,
where $\varphi:\,\cA\to\Hilb(\bP(\wedge^r(k^n))$
is a map induced by \eqref{KJYFKJTF}.
\end{Proposition}

\begin{proof}
As Lafforgue pointed out to us, this follows
from \cite[4.4,4.22]{Lafforgue03}.
\end{proof}

\begin{Review}\textsc{Structure map.}\label{LKJGKHFF}
We have the composition
$$
\oo \subset \cA \rightarrow \cA/\cA_{\emptyset}
$$
(where the last map is the stack quotient), which Lafforgue
calls the {\it structure map}.
In particular this endows
$\oo$ with a stratification by locally closed subschemes,
$\oo_{\uP}$ (the restriction of the corresponding
toric stratum of $\cA$), parameterized by matroid
decompositions $\uP$ of $\Delta(r,n)$. The stratum for
the trivial decomposition, $\emptyset$ (meaning the only
polytope is $\Delta(r,n)$) is an open subset
$$
\oo_{\emptyset} = X(r,n)\subset \oo
$$
which Lafforgue calls the {\it main} stratum.
Lafforgue proves that $\oo$ is projective, and thus
gives a {\it compactification} of $X(r,n)$ --
italics as his space is in general reducible, as
we observe in \eqref{sean3.7}.
\end{Review}

\begin{Review}\label{KLJFKHGDKHD}
We denote the closure of $\oo_\emptyset$ in $\oo$ by $\oX_L(r,n)$. There are
immersions
$$
X(r,n)\subset \oX_L(r,n) \subset \cA \subset \bP{//}_nH
$$
(the first and last open, the middle one closed)
and
$$
X(r,n)\subset \oX(r,n) \subset \bP// H
$$
(open followed by closed).
It follows that there exists a finite birational map
\begin{equation}\label{kjgfjhgfdjg}
\oX_L(r,n)\to \oX(r,n).\cooltag
\end{equation}
In particular, $\oX(r,n)$ and $\oX_L(r,n)$ have the same normalization
that we denote by $\tX(r,n)$.
\end{Review}

\begin{Review}\textsc{Toric family.}\label{LIUYFKJYF}
We denote the pullback of
$\tcA \rightarrow \cA$ to $\oo$ by
$\cT \rightarrow \oo$ ($\cT$~to denote toric).
By definitions, $\cT\subset\oo\times G(r,n)$.

Kapranov \cite[1.5.2]{Kapranov93} shows that $\oX(r,n)$
is isomorphic to the Hilbert quotient $G(r,n)///H$
and the natural Chow family
$$
\cT \rightarrow \oX(r,n),\quad \cT \subset \oX(r,n)\times G(r,n).
$$
is flat. The family
$\cT \rightarrow \oX_L(r,n)$ is the pullback of
$\cT \rightarrow \oX(r,n)$ along~\eqref{kjgfjhgfdjg}.

Let
$\cB,\cB_i \subset \cT$ be the restrictions of the boundary
divisors $\tcB,\tcB_i \subset \tcA$.
\end{Review}

\begin{Review}\label{sean2.6}\textsc{Family of visible contours.}
Let $G_e(r-1,n-1) \subset G(r,n)$ be the subspace of $r$-planes
containing the fixed vector $e=(1,\dots,1)$. Kapranov defines
the {\it family of visible contours}
$$
\cS = \cT \cap \bigl(\oX(r,n)\times G_e(r-1,n-1)\bigr)
\subset \oX(r,n)\times G(r,n).
$$
Kapranov shows that the family $\cS$ is flat, and that the
associated map
\begin{equation}\label{sean1.12}
\oX(r,n)\rightarrow \Hilb(G_e(r-1,n-1)) \cooltag
\end{equation}
is a closed embedding.
\end{Review}

There is a similar family over $\oo$
(Lafforgue calls it $\tbP(\cE)$):

\begin{Definition}\label{sean2.7}
Let $\cS \subset \cT$ be
the scheme theoretic intersection
$$
\cS := \cT \cap \lbrack \oo \times G_e(r-1,n-1) \rbrack \subset
\lbrack \oo \times G(r,n) \rbrack.
$$
\end{Definition}

$H$ acts on $
\tcA$, trivially on
$\cA$ and $\tcA \rightarrow \cA$ is $H$ equivariant. Thus
$H$ acts on $\cT$
(and trivially on $\oo$) so that $\cT \rightarrow \oo$ is
equivariant.

Let $\cB, \cB_i \subset \cS$ indicate the restriction of
$\cB, \cB_i \subset \cT$. We note $\cB \subset \cS$ is the
union of $\cB_i$, as the $n$ components of $\cB \subset \cT$
corresponding to the faces $x_i =0$ of $\Delta(r,n)$ are
easily seen to be disjoint from $G_e(r-1,n-1)$.

The fibres of $(\cS,\cB) \rightarrow \oo$ have singularities like (or
better) than those of $(\cT,\cB)$, as follows from the following
transversality result:

\begin{Proposition}[{\cite[pg xv]{Lafforgue03}}]\label{sean2.8}
The natural map
$$
\cS \rightarrow \cT/H
$$
to the quotient stack (or equivalently,
$\cS \times H \rightarrow \cT$) is smooth.
\end{Proposition}

\begin{proof}
We recall for the readers convenience Lafforgue's elegant
construction:
Let
$$
\cE \subset G(r,n) \times \aff^n
$$
be the universal rank $r$ subbundle, and let
$\czE \subset \cE$ be the inverse image under the
second projection of the open subset $H\subset \aff^n$
(i.e. the subset with all coordinates non-zero). $H$
obviously acts freely on $\czE$ and the quotient is canonically
identified with $G_e(r-1,n-1)$. This gives a smooth map
$$
G_e(r-1,n-1) = \czE/H\rightarrow G(r,n)/H
$$

Now for any $H$ equivariant $\cT \rightarrow G(r,n)$
the construction pulls back, yielding \eqref{sean2.8}.
\end{proof}

Note in particular that this shows

\begin{Corollary}\label{sean2.9}
$\cS \subset \cT$ is
regularly embedded, with normal bundle the pullback
of the universal quotient bundle of $G_e(r-1,n-1)$.
\end{Corollary}

\begin{Review}\label{sean2.10}\textsc{Fibers of $\cS$.}
A precise description of the fibres of $\cS$ is given
in \cite[Chapter 5]{Lafforgue03}. Here we recall a few points:

Let $S \subset T$ be a closed fibre of $\cS \subset \cT$ over
a point of $\oo_{\uP}$. We have
by the above a smooth {\it structure map} $S \rightarrow T/H$,
and so $S$ inherits a stratification from the orbit stratification
of $T/H$, parameterized by $P \in \uP$. In particular the
facets (maximal dimensional polytopes) of $\uP$ correspond to
irreducible components, and the stratum $S_P$ (which are
the points of $S$ that lie only on the irreducible component
corresponding to $P$) is the complement in $\bP^{r-1}$ to a
GIT stable arrangement (see~\eqref{sean1.15})
of $n$ hyperplanes with associated matroid polytope $P$ (see~\eqref{KJHGKG}).
The irreducible component itself is the log canonical
compactification of $S_P$, as follows for example from \eqref{sean2.12} below.
For $r=3$ this compactification is smooth, and described by \eqref{sean7.2}.
\end{Review}

\begin{Review}\label{sean2.11}\textsc{Log structures and toric stacks.}
For basic properties of log structures and toric stacks we refer to
\cite[\S5]{Olsson03}.  Any log structure
we use in this paper will be toric, i.e. the space will
come with an evident map to a toric variety and we endow
the space with the pullback of the toric log structure on
the toric variety. In fact, we do not make any use of
the log structure itself, only the bundles of log (and
relative log) differentials, all of which will be computed
by the following basic operation (our notation is chosen
with an eye to its immediate application):

Let $q:\tcA \rightarrow \cA$ be a map
of toric varieties so that the map of underlying tori is a
surjective homomorphism, with kernel $H$. We have the
smooth map
$$
\tcA \rightarrow \tcA/H
$$
(where the target is the stack quotient), and in particular
its relative cotangent bundle, which is canonically identified
with a trivial bundle with fibre the dual of the Lie algebra
to $H$. We denote the bundle
\begin{equation}\label{lkhgkhjfgkf}
\Omega^1_q(\log) = \Omega^1_{\tcA/(\tcA/H)}
\cooltag
\end{equation}
as $q$ is log smooth and this is its
bundle of relative log differentials, as follows from
\cite[5.14]{Olsson03} and \cite[3.7]{Olsson}.

For a map
$\oo \rightarrow \cA$, consider the pullback
$$
\cT := \tcA \times_{\cA} \oo \rightarrow \oo
$$
Then \eqref{lkhgkhjfgkf} pulls back to
the relative cotangent bundle for
$$
\cT \rightarrow \cT/H.
$$
$\cT \rightarrow \oo$ is again log smooth, with this
(trivial) bundle of relative log differentials.

Now suppose $\cS \subset \cT$ is a closed subscheme, so that
the map $\cS \times H \rightarrow \cT$, or
equivalently, $\cS \rightarrow \cT/H$, is smooth.
Then the relative cotangent bundle for
$$
\cS \rightarrow \cT/H
$$
is a quotient of the pullback of $\Omega^1_{\cT/(\cT/H)}$,\quad
$p:\, \cS \rightarrow \oo$ is log smooth, with bundle of
relative log differentials
$$
\Omega^1_p(\log) = \Omega^1_{\cS/(\cT/H)}.
$$
\end{Review}

\begin{Theorem}\label{sean2.12}
The visible contour family
$p:\cS \rightarrow \oo$ is log
smooth. Its bundle of log differentials
$$
\Omega^1_p(\log) = \Omega^1_{\cS/(\cT/H)}
$$
is a quotient of the pullback of
$\Omega^1_{\tcA/(\tcA/H)}$,
which is the trivial bundle $\tcA \times V_n$. Fibres
$(S,B)$ are semi-log canonical, and the restriction of
the Pl\"ucker polarisation to $S \subset G_e(r-1,n-1)$
is $\cO(K_S + B)$.
\end{Theorem}

\begin{proof} Let $(S,B) \subset (T,B)$ be closed fibres
of $(\cS,\cB) \subset (\cT,\cB)$. $(T,B)$ is semi-log
canonical, and $\cO(K_T + B)$ is canonically trivial,
e.g. by \cite[3.1]{Alexeev}.
$(S,B)$ is now semi-log canonical by \eqref{sean2.8}, and by adjunction
$\cO(K_S + B)$ is the determinant of its normal bundle,
which is the Pl\"ucker polarisation by \eqref{sean2.9}. The
other claims are immediate from \eqref{sean2.8} and the general
discussion \eqref{sean2.11}.
\end{proof}

The initial motivation for this paper was
the elementary observation:

\begin{Proposition}\label{sean1.8}
$X(r,n)$ is minimal of log general type.
\end{Proposition}

\begin{proof}
We show that the first log canonical
map on $X(r,n)$ is a regular immersion.
Fixing the first $r+1$ hyperplanes identifies
$X(r,n)$ with an open subset of
$U^{n-(r+1)}$, where $U \subset \bP^{r-1}$ is the
complement to $B$, the union of $r+1$ fixed hyperplanes
in linear general position.
$K_{\bP^{r-1}} + B = \cO(1)$, so the first
log canonical map on $U$ is just the inclusion
$U \subset \bP^{r-1}$, in particular an immersion. The
result follows easily.
\end{proof}

We have the following criterion to guarantee
$$X(r,n)\subset \oX_L(r,n)$$
is a log minimal model.
Let $T_p(\log)$ be the dual
bundle to $\Omega^1_p(\log)$ on $\cS$ -- i.e. the relative
tangent bundle to $\cS \rightarrow \cT/H$.

\begin{Theorem}\label{sean2.16}
If $R^2p_*(T_p(\log))$ vanishes at a point of
$\oX_L(r,n)\subset \oo$, then
$\oo \rightarrow \cA/\cA_\emptyset$ is smooth,
$\oX_L(r,n) = \oo$, $\oo$ is normal, and the pair
$(\oX_L(r,n),B)$ has toroidal singularities, near the
point.

If $R^2p_*(T_p(\log))$ vanishes identically along
$\oX_L(r,n)$, then the sheaf
$$\Omega^1_{\oX_L(r,n)/k}(\log B)$$ {\rm(}defined
in \eqref{sean9.6}{\rm)}
is locally free, globally generated, and its determinant,
$\cO(K_{\oX_L} + B)$, is globally generated and big.
In particular
$$
X(r,n) \subset \oX_L(r,n)
$$
is a log minimal model.
\end{Theorem}

\begin{proof}
By \cite[4.25.ii,5.15]{Lafforgue03},
vanishing of $R^2$ implies the structure map
is smooth. Now suppose the structure
map is smooth along $\oX_L(r,n)$. The bundle of log differentials
for the toric log structure on a normal toric variety is
precisely the bundle \eqref{sean9.6}, which implies the analogous statement
for $(\oX_L(r,n),B)$.
The bundle of
differentials is the cotangent bundle of the structure
sheaf, and thus a quotient of the cotangent bundle to
$$
\tcA \rightarrow \tcA/\tcA_{\emptyset}
$$
which by \eqref{sean2.11} is canonically trivial, whence the global
generation. Now $K_{\oX_L(r,n)} + B$ is big by \eqref{sean1.8}.
\end{proof}

\begin{Remark}
If the conditions of the theorem hold, then to show
$\oX_L(r,n)$ is the log canonical model, it remains to show
$K_{\oX_L} + B$ is ample, not just big and nef. We have
proven this for $\oX_L(3,6)$, by restricting to the boundary.
We expect it will hold whenever \eqref{sean2.16} applies, i.e. in
the cases of \eqref{sean1.6}.
\end{Remark}


\begin{Remark}
When vanishing holds in \eqref{sean2.16} we have generating global
sections of $\Omega^1_{\oX_L(r,n)/k}(\log B)$ which give a map
$$
\oX_L(r,n) \rightarrow G\left((r-1)(n-r-1),\binom{n}{r} - n\right).
$$
For $r=2$ the sections give a basis of global sections, and
we have checked the map is a closed embedding. It would be
interesting to know the defining equations. In this case the
log canonical series is very ample, and the embedding factors
through this embedding into the Grassmannian. The log canonical
embedding itself is quite nice, see~\eqref{sean2.17}.
\end{Remark}

\begin{Review}\textsc{When is $\Psi$ bijective?}
Here we resume a notation of~\eqref{HGPOU} and give a technical condition
that implies $\Psi|_{(\bP{//}_nH)_\uP}$ is bijective for a polyhedral decomposition $\uP$.
Until the end of this section, we assume that $\uP$ is unimodular.

The construction is a variation of the Ishida's complex of $\bZ$-modules,
see~\cite{Oda}. Let $\uP^i$ be the set of $i$-codimensional faces of polytopes in $\uP$
that do not belong to the boundary $\partial P$. We fix some orientation of each $Q\in\uP^i$.
Let $A$ be an abelian group.
Consider the complex $C^\bullet_{\Aff}(\uP,A)$ with
$C^i_{\Aff}(\uP,A)=\oplus_{Q\in\uP^i}\Aff(Q,A)$,
where $\Aff(Q,A)$ is the group of affine maps $Q\to A$. The differential $d^i:\,C^i_{\Aff}(\uP,A)\to C^{i+1}_{\Aff}(\uP,A)$
is a direct sum of differentials $d^{Q,R}$ for $Q\in\uP^i$, $R\in\uP^{i+1}$.
If $R$ is not a face of $Q$ then $d^{Q,R}=0$. Otherwise, $d^{Q,R}$ is the restriction
map $\Aff(Q,A)\to\Aff(R,A)$ taken with a negative sign if the fixed orientation
of $R$ is opposite to the orientation induced from $Q$.
Let $H^\bullet_{\Aff}(\uP,A)$ be the cohomology of $C^\bullet_{\Aff}(\uP,A)$.
It is clear that $H^0_{\Aff}(\uP,A)$ is the set of piecewise affine functions $P\to A$.
\end{Review}

\begin{Proposition}\label{JKFKJGF}
If $H^1_{\Aff}(\uP,\bZ)=0$ then $\Psi|_{(\bP{//}_nH)_\uP}$ is bijective.
\end{Proposition}

\begin{proof}
We identify $\cH$ with maps $P\to\bG_m$. Elements of $\cH$ of order $N$ are maps $P\to\mu_N$
and any map $a:\,P\to\bZ$ gives a 1-PS $z\mapsto \{p\mapsto z^{a(p)}\}$.

Let $X\in(\bP//H)_\uP$ be as in \eqref{brokentoric}. Let $x\in\Psi^{-1}(X)$.
We claim that $\cH x \to \cH X$ is bijective.
Since $\cH_x\subset \cH_X$,  it sufices to prove that the stabilizer $\cH_X$ is connected.

Let $h\in\cH_X$. Then
$h\in\cH_{X_{P'}}$ for any $P'\in\uP$.
But if $e$ is a generic point of $X_{P'}$ then
$$\cH_{X_{P'}}=\{h\in\cH\,|\,h\cdot e\in X_{P'}\}=\{h\in\cH\,|\,\exists h_{P'}\in H,\ h\cdot e=h_{P'}\cdot e\}.$$
It follows that
$h(p)=h_{P'}(p)$ and hence $h$ is affine on each $P'$.
We see that $\cH_X=H^0_{\Aff}(\uP,\bG_m)$.

It is enough to show that any element $h\in\cH_X$
of a finite order $N$  embeds in a 1-PS $\gamma\subset\cH_X$.
So let $h\in H^0_{\Aff}(\uP,\mu_N)$.
We have the exact sequence $$0\to C^\bullet_{\Aff}(\uP,\bZ)\mathop{\longrightarrow}^{\cdot N} C^\bullet_{\Aff}(\uP,\bZ)\to C^\bullet_{\Aff}(\uP,\mu_N)\to0.$$
Since $H^1_{\Aff}(\uP,\bZ)=0$,
there exists an element of $H^0_{\Aff}(\uP,\bZ)$ that maps to $h$.
The corresponding 1-PS contains $h$ and belongs to the stabilizer $\cH_X$.
\end{proof}

\begin{Definition}\label{KLJGKF}
A decomposition $\uP$ is called {\em central} if $\uP^0=\{C,S_1,\ldots,S_r\}$,
where $S_i\cap S_j\subset\partial P$. We call $C$ the central polytope.
Let $U_{\uP}\subset \bP{//}_nH$ be an affine open toric subset with fan $C(\uP)$. It contains
$(\bP{//}_nH)_\uP$ as the only closed orbit.
\end{Definition}

\begin{Corollary}\label{sean3.4}
If $\uP$ is unimodular and central then $\Psi(U_\uP)$ is quasi-smooth, i.e.~$\Psi|_{U_{\uP}}$ is bijective
and $U_{\uP}$ is smooth.
\end{Corollary}

\begin{proof}
To show that $\Psi|_{U_{\uP}}$ is bijective, it suffices to prove that
$\Psi|_{(\bP{//}_nH)_\uP}$ is bijective. Indeed, other strata in
$U_{\uP}$ correspond to decompositions coarser than $\uP$, which are automatically
unimodular and central, so we can use the same argument.
It is clear that $\uP^{1}=\{F_1,\ldots,F_p\}$ is the set of codimension~$1$ faces of~$C$ that are not on the boundary of $P$. We want to use~\eqref{JKFKJGF}.
Let $c\in C^1_{\Aff}(\uP,\bZ)$, so $c=(f_1,\ldots,f_p)$, where $f_i$ is affine on $F_i$.
For each $i$, we have $F_i=C\cap S_j$ for some $j$.
We can choose $g_j\in\Aff(S_j,\bZ)$ which restricts on $f_i$
(taking into account the orientation). Then $c$ is equal to the differential of the cochain~$\tilde c$,
where $\tilde c(C)=0$ and $\tilde c(S_i)=g_i$.

For the second statement, we have to show that $C(\uP)=(\bZ_{\ge0})^r$ upto global affine functions.
Let $f\in C(\uP)$. Then $f$ is a concave locally affine function. So $f-f|_C$ is a concave locally affine
function that vanishes on $C$. Let $f_i$, $1\le i\le r$, be a primitive (i.e.~not divisible by an integer)
concave locally affine function that vanishes on $P\setminus S_i$.
Then $f-f|_C$ is a linear combination of $f_i$'s with non-negative coefficients.
\end{proof}

\section{Singularities of $(\oX(r,n), B)$}\label{sean3}
In this section we prove (1.13). The very simple idea is as follows:
The notion of log canonical pair $(\oX,\sum B_i)$ generalises normal
crossing. In particular, if all the irreducible components $B_i$
are $\bQ$-Cartier, then log canonical implies at least that
the intersection of the $B_i$ has the expected codimension, see 
\eqref{sean3.13} below. We prove (1.13) by observing that well
known configurations 
give points of $\oX(r,n)$ lying on too many boundary divisors. The
main work is to show that these points are actually in the 
closure of the generic stratum, and that the boundary divisors
are Cartier near these points. 

\begin{Review}\textsc{Matroid polytopes \cite{GMS87}.}\label{KJHGKG}
Let $N=\{1,\ldots,n\}$. A {\em matroid} $C$ is $N$
together with a nonempty family of {\em independent subsets} of $N$ such that
any subset of an independent subset is independent and
all maximal independent subsets contained in $I$ have the same number of elements
for any $I\subset N$.
Maximal independent subsets of $N$ are called {\em bases} of the matroid.
Obviously, a subset of a matroid is a matroid
(with induced collection of independent subsets).
The {\em rank} of a matroid is equal to the number of elements in any base.
A matroid $C$ of rank $r$ gives rise to a {\em matroid polytope}
$P_C\subset \Delta(r,n)$,  a convex hull of vertices $e_{i_1}+\ldots+e_{i_r}$ for any base
$\{i_1,\ldots,i_r\}\subset N$.
For a subset $I \subset N$, we write
$x_I=\sum_{i \in I} x_i$.
We consider $x_I$ as a function on $\Delta(r,n)$, in particular
$x_I=r-x_{I^c}$. It is known that
$P_C$ is defined by inequalities $x_I\le\rank I$.
\end{Review}

\begin{Review}\textsc{Realizable matroids.}\label{sean3.1}
Here is the main example: Let $C=\{L_i\}_{i\in N}$ be a configuration of $n$ hyperplanes in $\bP^{r-1}$.
Then independent subsets of the corresponding {\em realizable matroid} (denoted by the same letter $C$)
are subsets of linearly independent hyperplanes.
$C$ has rank $r$ if there is at least one independent $r$-tuple.
Let $X_C(r,n)$ be the corresponding moduli space, i.e. $N$-tuples
of hyperplanes with incidence as specified by~$C$ modulo $\PGL_r$.
The corresponding matroid polytope $P_C$ has a maximal dimension iff $C$ is GIT stable.
By the {\em multiplicity} of a point $p \in \bP^{r-1}$ with respect to~$C$ we mean the
number of hyperplanes in $C$ that contain $p$, i.e.~the usual geometric
$\mult_pC$ if we view $C$ as a divisor in  $\bP^{r-1}$.
\end{Review}

\begin{Review}\textsc{Gelfand--Macpherson correspondence.}\label{KHGKHGJ}
Let $C$ be as in \eqref{sean3.1}. Consider an $r\times n$ matrix $M_C$ with columns
given by equations of hyperplanes of $C$ (defined upto a scalar multiple).
The row space of $M_C$ gives a point of $G(r,n)$.
Thus $X_C(r,n)$ is identified with
the quotient of (the reduction of) a thin Schubert cell
$G^{P_C, 0}(r,n)/H$ (see \eqref{KJHGKHGJG}).
So we see that for any $x\in G(r,n)$, $\Supp(x)\subset\Delta(r,n)$
is a matroid polytope of a realizable matroid, in particular the Lafforgue's stratum $\oo_\uP$ \eqref{LKJGKHFF}
is empty if a matroid decomposition $\uP$ contains non-realizable matroids.
\end{Review}

\begin{Review}\textsc{Divisor $B_I$.}
It is easy to see that $\{x_I\le k\}$ is a matroid polytope for any $0<k<r$.
The corresponding configuration is as follows:
the only condition we impose is
$$\codim\bigcap_{i\in I}L_i=k.$$
This polytope has full dimension iff $|I|>k$

It follows that if $|I|>k$ and $|I^c|>r-k$
then there is a matroid decomposition of $\Delta(r,n)$
with two polytopes $\{x_I \ge k\}$ and $\{x_{I}\le k\}$.
The corresponding stratum of $\cA$ is maximal among boundary
strata. 
We denote its closure (and corresponding subschemes
of $\oo$, $\oX_L(r,n)$, etc.) by~$B_I$.

An example is shown in \eqref{sean1.19}, where $r=3$, $n=5$, $k=1$, $I=\{2,3,4\}$, $I^c=\{1,5\}$.
In the configuration with polytope $\{x_{I^c} \le 1\}$, lines of $I^c$ are identified and lines of $I$ are generic.
In the configuration $\{x_{I}\le 2\}$, the lines of $I$ have a common point
of incidence, and lines of $I^c$ are generic.
\end{Review}

\begin{Review}\textsc{Central configurations and matroids.}
Let $\cI$ be an index set, and
for each $\alpha\in\cI$, $I_\alpha\subset N$ a subset, such that
$|I_\alpha|\ge r$  and
\begin{equation}\label{:LKJHGLKHG}
|I_\alpha \cap I_\beta| \leq r-2\quad\hbox{\rm for}\quad \alpha\neq\beta.\cooltag
\end{equation}
Let's call $S\subset N$ independent if $|S|<r$ or $|S|=r$ and $S\not\subset I_\alpha$ for any $\alpha\in\cI$.
\end{Review}

\begin{Proposition}
This gives a structure of a matroid on $N$.
\end{Proposition}

\begin{proof}
We only have to check that for any $S\subset N$, all maximal independent subsets in $S$
have the same number of elements. It suffices to prove that if $S$ contains an independent set $T$, $|T|=r$,
then any independent subset $R\subset S$ can be embedded in an independent subset with $r$ elements.
We can assume that $|R|=r-1$. If $R\not\subset I_\alpha$ for any $\alpha\in\cI$, then we can just add any element to $R$.
If $R\subset I_\alpha$ for some $\alpha$ then this $\alpha$ is unique by~\eqref{:LKJHGLKHG}
and we add to $R$ an element of $T$ that is not contained in $I_\alpha$.
\end{proof}

We call matroids of this form {\em central}.

\begin{Definition}\label{:LKJGLKH}\label{sean3.2}
We say that $C$ is a {\em central configuration} if a pair $(\bP^{r-1},C)$ has normal crossings
on the complement to a $0$-dimensional set. If $r=3$, it simply means that there are no double lines.
Let $\cI\subset\bP^{r-1}$ be the set of points of multiplicity at least~$r$.
Then a matroid of $C$ is a central matroid that corresponds to subsets
$I_\alpha \subset N$ of hyperplanes containing $\alpha\in \cI$.

A polytope $P_C$ of a central matroid $C$ is given by inequalities
$x_{I_\alpha} \leq r-1$ for all $\alpha\in\cI$.
Let $P_\alpha \subset \Delta(r,n)$ be the matroid polytope
$x_{I_\alpha}\geq r-1$.
Let $\ucI=\{P_C,P_\alpha\}_{\alpha\in\cI}$.
\end{Definition}

\begin{Lemma}\label{sean3.3}
$\ucI$ is a central decomposition of $\Delta(r,n)$ {\rm(}see~\eqref{KLJGKF}
for the definition{\rm)}
with central polytope $P_C$.
For each subset
$\cI'\subset\cI$, $\ucI'$ is a matroid decomposition, coarser
than $\ucI$ and all matroid decompositions coarser than
$\ucI$ occur in this way.
\end{Lemma}

\begin{proof}
To show that $\ucI$ is a central decomposition,
it suffices to check that $P_\alpha\cap P_\beta$ is on the boundary for any $\alpha\neq\beta$
(this will imply, in particular, that any interior point of any wall $\{x_{I_\alpha}=r-1\}$
belongs to exactly two polytopes, $P_C$ and $P_\alpha$).
Assume that $x\in P_\alpha\cap P_\beta$ is an interior point of $\Delta(r,n)$.
Then $x_{I_\alpha\cap I_\beta}<r-2$ by~\eqref{:LKJHGLKHG}
(otherwise $x_i=1$ for any $i\in I_\alpha\cap I_\beta$ and therefore $x$ is on the boundary).
Therefore, $x_{I_\alpha\setminus I_\beta}=x_{I_\alpha}-x_{I_\alpha\cap I_\beta}>1$ and
$x_{I_\alpha\cap I_\beta}= x_{I_\alpha\setminus I_\beta}+x_{I_\beta}>r$. Contradiction.

Any matroid decomposition coarser than $\ucI$ is obviosly central and can be obtained
by combining $P_C$ with several $P_{\alpha}$'s. This has the same effect as
taking these $\alpha$'s out of $\cI$.
\end{proof}

\begin{Proposition}\label{<JGKJGFGJF}
Let $U_{\ucI}$ be the affine open toric subset of $\cA$ as in~\eqref{KLJGKF}.
Then $U_{\ucI}$ is smooth and bijective to 
$\Psi(U_{\ucI})\subset\bP(\Lambda^r k^n)//H$. Let $U=U_{\ucI}\cap\oo$.
$U\subset \oX_L(r,n)$  maps finitely and homeomorphically onto 
its image in 
$\oX(r,n)$.
\end{Proposition}

\begin{proof} Follows from~\eqref{sean3.4}.
\end{proof}

\begin{Remark}
It is possible that $\ucI$ corresponds to a central
configuration~$C$ but central polytopes of decompositions coarser than $\ucI$ are not realizable.
In other words, some multiple points of $C$ may be forced by the set of
other multiple points (the reader may wish to examine picture~\eqref{sean3.7} from this perspective). All configurations used in the proof of \eqref{sean3.18} are of this sort.
\end{Remark}

\begin{Definition}\label{sean3.5}
Fix a hyperplane $L \subset \bP^{r-1}$. For each
subset $J \subset N$, $|J| \geq r$, Let $Q_J$ be the moduli
space of $J$-tuples of hyperplanes, $L_j$, $j \in J$ in $\bP^{r-1}$ such
that the entire collection of hyperplanes, together with $L$ is
in linear general position, modulo automorphism of $\bP^{r-1}$ preserving
$L$.
\end{Definition}

Note $Q_J$ is a smooth variety, of dimension $(r-1)(|J|-r)$. Intersecting
with the fixed hyperplane $L$ gives a natural smooth surjection
$$
Q_J \rightarrow X(r-1,|J|).
$$

\begin{Lemma}\label{sean3.6}
Let $C$, $\cI$, $\ucI$ be as in \eqref{:LKJGLKH}. For each
$\alpha\in\cI$ we have a natural map $X_C(r,n) \rightarrow X(r-1,|I_\alpha|)$,
taking the hyperplanes through $a$. Let
$$M = \prod_{\alpha \in \cI}X(r-1,|I_\alpha|)\quad \text{and}\quad Q = \prod_{\alpha \in \cI} Q_{I_\alpha}.$$
There is a natural identification
$$
\oo_{\ucI} = X_C(r,n) \times_M  Q.
$$
In particular
$$
\dim(\oo_{\ucI}) = \dim(X_C(r,n)) + \sum_{\alpha \in \cI}(|I_a| - r).
$$
\end{Lemma}

\begin{proof} This is immediate from \cite[\S 3.6]{Lafforgue03}.
\end{proof}

Next we demonstrate that Lafforgue's space $\oo$ is reducible:

\begin{Proposition}\label{sean3.7}
Let $C$ be the following
configuration of $6m-2$ lines in $\bR^2$:

\bigskip
\epsfxsize=0.5\textwidth
\centerline{\epsfbox{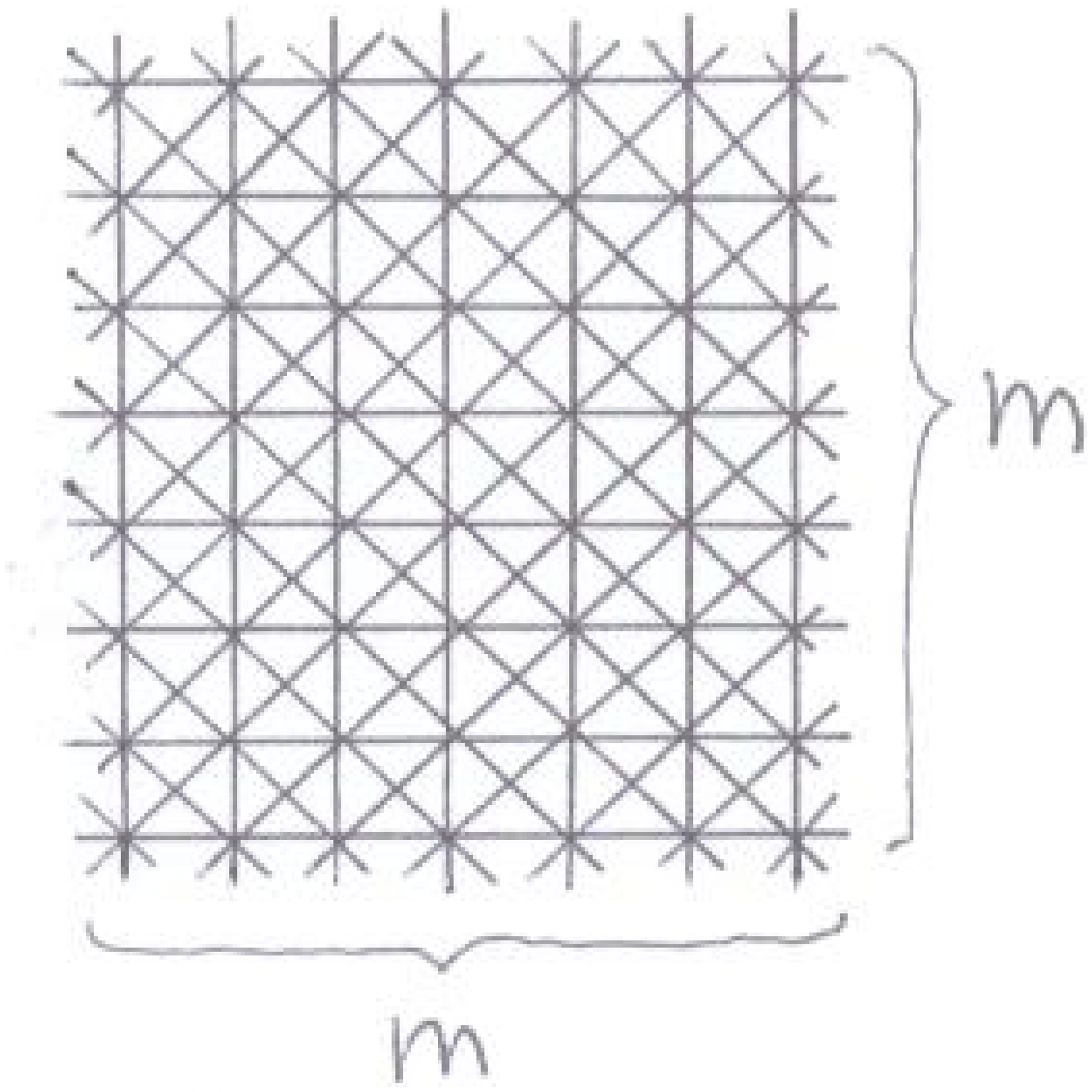}}
\bigskip

\noindent
Let $\cI$ be its multiple points, as in \eqref{:LKJGLKH}. Then
$$
\dim\left(\oo^{\Delta(3,6m-2)}_{\ucI}\right) \geq m^2
$$
and
$\oo^{\Delta(3,n)}$ is not irreducible for large $n$.
\end{Proposition}

\begin{proof} The configuration $C$ has at least $m^2$ points of multiplicity $4$,
so the inequality is immediate from \eqref{sean3.6}. The final remark
follows as the main component $\oX_L(3,6m-2)$ of $\oo^{\Delta(3,6m-2)}$ has
dimension $12m -12$.
\end{proof}

However, for a large class of central configurations, the stratum $\oo_{\ucI}$
belongs to the closure of the main stratum:

\begin{Review}\label{sean3.8}\textsc{Lax configurations.}
We say that a central configuration $C$ is {\em lax} if there is
a total ordering on $N$ so that for each $i \in N$,
points on $L_i$ of multiplicity {\em greater} than $r$ with respect to $N_{\le i}$
are linearly independent.
For example, a configuration in \eqref{sean3.7} is not lax for $m\ge4$.
\end{Review}

\begin{Theorem}\label{sean3.12}
Notation as in \eqref{sean3.2}. Assume $C$ is lax.
\begin{enumerate}
\renewcommand\theenumi{\arabic{section}.\arabic{Theorem}.\arabic{enumi}}
\item\label{KJFKJFKKHGD}
The stratum
$\oo_{\uI}$ is contained (set theoretically)
in $\oX_L(r,n)\subset \oo$.

\item\label{:LKJH:KJH}
Let $U=U_{\ucI}\cap\oo$, where $U_{\ucI}\subset\cA$ is the smooth
toric affine open set of~\eqref{<JGKJGFGJF}. Let
$\tU \rightarrow U$ be the normalisation. 
Then $U$ is an irreducible open factorial subset of
$\oX_L(r,n)\subset \oo$, smooth in codimension one.
Moreover  the boundary strata $B_{I_a}$ are Cartier, generically
smooth, 
and irreducible on $U$, their union is the boundary, and
their scheme-theoretic intersection is the stratum 
$\oo_{\ucI}$.

\item\label{LKHGKKG}
Let $\tU \rightarrow U$ be the normalisation, and $\tB \subset \tU$
the reduction of the inverse image of $B$. 
If $K_{\tU} + \tB$ is log canonical at a point in the
inverse image of $p \in \oo_{\ucI}$ then
the stratum has pure codimension $|\cI|$ in $U$ near $p$, i.e.
$$
\sum_{\alpha \in \cI} (|I_\alpha| - r+1) + \dim{X_C(r,n)}=n(r-1)-r^2+1
$$
near $p$.
\end{enumerate}
\end{Theorem}

We postpone the proof until the end of this section.
First we show that
$\tX(3,n)$ with its boundary fails to be
log canonical for $n \geq 9$
(for $n \geq 7$ in characteristic $2$) and that
$\tX(4,n)$ is not log canonical for $n\ge8$.

\begin{proof}[Proof of~\eqref{sean1.5}]
Consider the Brianchon--Pascal configuration \cite{HC,Do}
of 9 lines with $|\cI|=9$ and $|I_\alpha| = 3$ for all $\alpha$:

\bigskip
\epsfxsize=0.5\textwidth
\centerline{\epsfbox{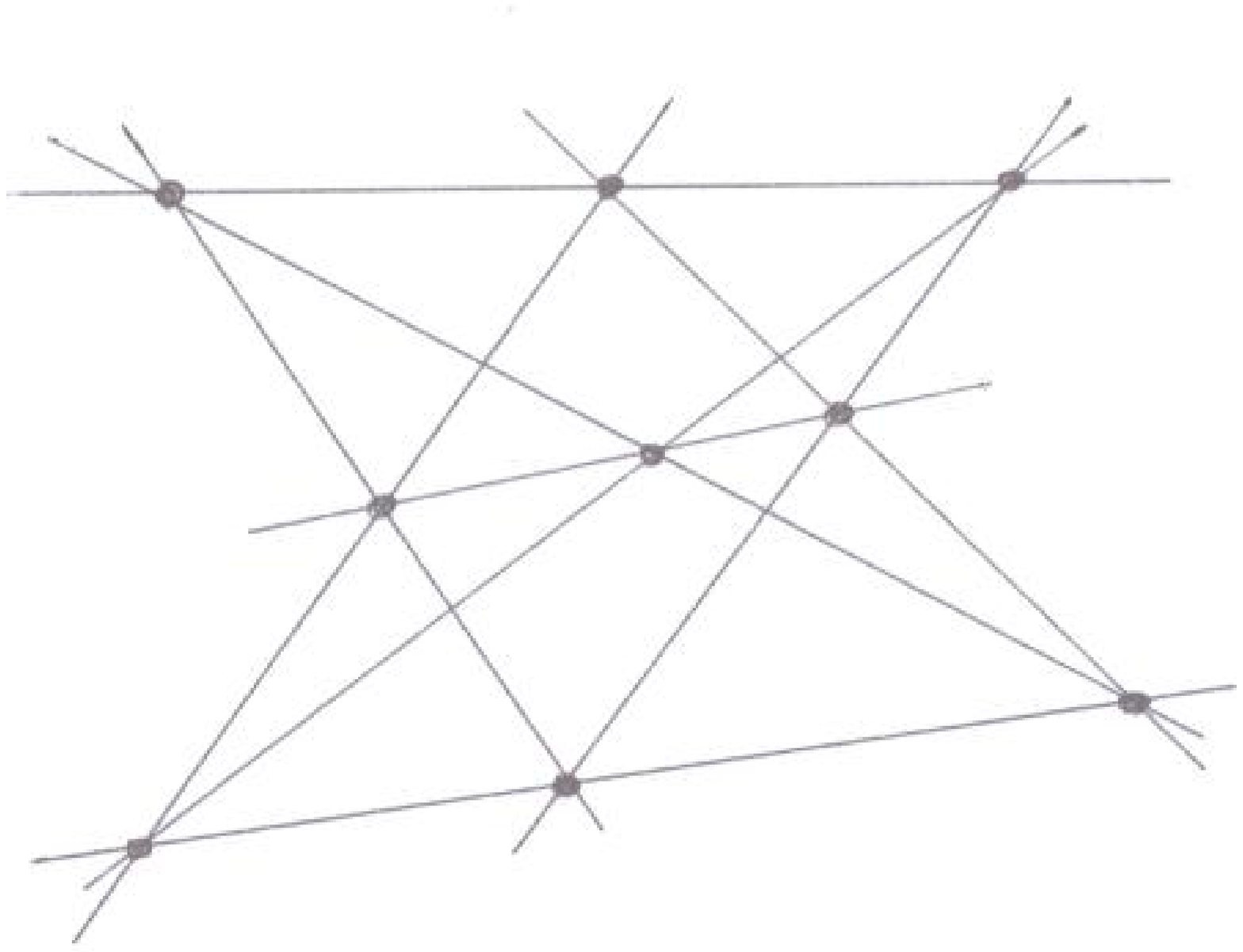}}
\bigskip
It is easy to compute that $\dim X_C(3,9)=2$.
Now apply \eqref{sean3.12}: the LHS in \eqref{LKHGKKG} is equal to $11$ but the RHS is $10$.
If $n \geq 10$ add generic lines.

There is an even better configuration of $9$ lines with
$|\cI| = 12$ and $|I_\alpha| = 3$ for all $\alpha$. It can be obtained
as follows: Fix a smooth plane cubic. Every line containing
two distinct inflection points contains exactly three. This
gives a configuration of $12$ lines. Furthermore each
inflection point lies on exactly $3$ lines, and these
are all the intersection points of the configuration.
This is the famous Hesse Wendepunkts-configuration \cite{HC,Do}. Let $C$ be the dual configuration.
Now apply \eqref{sean3.12}: the LHS in \eqref{LKHGKKG} is equal to $12$ but the RHS is $10$.
If $n \geq 10$ add generic lines.

For the characteristic two, use the Fano configuration,
\cite[4.5]{GMS87} and argue as above:
the LHS in \eqref{LKHGKKG} is equal to $7$ but the RHS is $6$.

In $(4,8)$ case, take the configuration of $8$ planes
in $\bP^3$ given by the faces of the octahedron.
There are $12$ points of multiplicity $4$ (i.e. lying
on $4$ of the planes), while $\tX(4,8)$ is $9$ dimensional.
If $n \geq 9$ add generic planes.
\end{proof}

\begin{Theorem}\label{sean3.17}
The boundary
strata of $(\oX(3,n),B)$ for lax configurations have arbitrary singularities,
i.e. their reductions give reductions of all
possible affine varieties defined over $\bZ$ (up to products with $\aff^1$).
\end{Theorem} 

\begin{proof}
By \eqref{sean3.12}, it suffices to prove that
$\oo^{\Delta(3,n)}_{\ucI}$ for
lax configurations $C,\ucI$, satisfy Mnev's theorem \cite[1.14]{Lafforgue03}. I.e.
given affine variety $Y$ over $\bZ$ there are integers $n,m$ and an
open set $U \subset Y \times \aff^m$, with $U \rightarrow Y$
surjective, and a lax configuration $C$ with $n$
lines such that $U$ is isomorphic to the reduction of the Lafforgue stratum
$\oo^{\Delta(3,n)}_{\ucI}$.
One can follow directly the proof of Mnev's theorem:
Lafforgue constructs an explicit configuration which encodes
the defining equations for $Y$, and it is easy to check this
configuration is lax.
The ordering of lines (in Lafforgue's notation) should be as follows:
lines $[0,1_\alpha,P_\alpha,\infty_\alpha]$ and the infinite line should go first
(at the end of the process there will be many points of multiplicity $>3$ along them),
then take all auxilliary lines in the order of their appearance in the Lafforgue's construction.
\end{proof}

Now we proceed with the proof of \eqref{sean3.12}.

\begin{Review}\textsc{Face maps and cross-ratios.}\label{KGLJGLKYF}
The collection of $X(r,n)$ has a hypersimplicial structure: there are obvious maps
$B_i:\,X(r,n)\to X(r,n-1)$ (dropping the $i$-th hyperplane)
and $A_i:\,X(r,n)\to X(r-1,n-1)$ (intersecting with the $i$-th hyperplane).
These maps extend to maps of Chow quotients \cite[1.6]{Kapranov93}, and to maps
of Lafforgue's varieties $\oX_L\subset\oo\subset\cA$, \cite[2.4]{Lafforgue03}.
For $\cA$, these maps are just restrictions of face maps \eqref{KJHGKHGJG}
corresponding to faces $\{x_i=0\}\simeq\Delta(r,n-1)$ and
$\{x_i=1\}\simeq\Delta(r-1,n-1)$.

In particular, let $V, W\subset N$ be subsets such that $|V|=4$, $|W|=r-2$, $V\cap W\ne\emptyset$.
Then dropping all hyperplanes not in $V\cup W$ and intersecting
with all hyperplanes in $W$ gives cross-ratio maps
\begin{equation}\label{KJGKGHG}
CR_{V,W}:\,X(r,n)\to X(2,4)=M_{0,4}=\bP^1\setminus\{0,1,\infty\}\cooltag
\end{equation}
and
$$
CR_{V,W}:\,\oo\to\oX(2,4)=\oM_{0,4}=\bP^1.$$
It follows that $CR_{V,W}(\oo_\uP)\subset\bP^1\setminus\{0,1,\infty\}$
if and only if $\uP$ does not break
$$\Delta_{V,W}(2,4)=\bigcap_{i\not\in V\cup W}\{x_i=0\}\bigcap\bigcap_{i\in W}\{x_i=1\}.$$
$\Delta(2,4)$ is an octahedron and values $\{0,1,\infty\}$ correspond to three decompositions
of $\Delta(2,4)$ into two pyramids.

To write \eqref{KJGKGHG} as a cross-ratio, let $V=\{i_1i_2i_3i_4\}$ with
$i_1<i_2<i_3<i_4$.
Let $L_1,\ldots,L_n$ be a collection of hyperplanes in $X(r,n)$.
Consider an $r\times n$ matrix $M$ with columns given by equations of these 
hyperplanes.
Then
$$CR_{V,W}(L_1,\ldots,L_n)={\Det_{i_1i_2W}\Det_{i_3i_4W}\over\Det_{i_1i_3W}\Det_{i_2i_4W}},$$
where each $\Det_T$ is an $r\times r$ minor of $M$ with columns given by $T$.
\end{Review}

Let $C=\{L_1,\ldots,L_n\}$ be any configuration as in~\eqref{sean3.1} and
let $x_0\in G(r,n)$ be a point that corresponds to $C$
under the Gelfand--Macpherson transform.
Let $(\cX,x_0)\subset G(r,n)$ be a pointed curve such that 
$\cX\cap G^0(r,n)\ne\emptyset$.
Let $F:\,G(r,n)^0\to X(r,n)$ be the canonical $H$-torsor.
Then
$$p_0=\lim\limits_{x\to x_0}F(x)\in X_L(r,n)$$
belongs to $\oo_\uP$, where $\uP$ is a
matroid decomposition of $\Delta(r,n)$ containing~$P_C$.
Indeed, it is clear that $x_0$ is contained in the fiber of the 
universal family~\eqref{sean2.2} over $p_0$,
so $P_C=\Supp(x_0)$ is in $\uP$.

\begin{Proposition}\label{JFJFFJFJF}
Let $C$ be central as in~\eqref{sean3.2}.
If $$\lim\limits_{x\to x_0}CR_{V,W}(x)\not\in\{0,1,\infty\}$$
for any $W\subset I_\alpha$, $|V\cap I_\alpha^c|=1$, $\alpha\in\cI$ then 
$\uP=\ucI$.
\end{Proposition}

\begin{proof}
Any decomposition containing $P_C$ is a refinement of~$\ucI$.
So it remains to prove the following combinatorial statement: any 
realizable matroid decomposition $\uP$
refining $\ucI$ is equal to $\ucI$ provided that 
$\uP\cap\Delta(2,4)=\Delta(2,4)$ for any face
$\Delta(2,4)\subset\Delta(r,n)$ that belongs to the boundary of some 
$P_\alpha$ and such that
exactly one face of this octahedron $\Delta(2,4)$ belongs to the wall 
$x_{I_\alpha}=r-1$
(this is a condition equivalent to $W\subset I_\alpha$, $|V\cap I_\alpha^c|=1$).

Restrictions of $\uP$ and $\ucI$ to the faces of $\Delta(r,n)$ have the same form.
Also, if $r=2$ then the claim follows, for example, from the
explicit description of matroid decompositions of 
$\Delta(2,n)$~\cite[1.3]{Kapranov93},
so we can argue by induction and it remains to prove the following:
any realizable matroid decomposition $\uP$ refining $\ucI$ is equal to 
$\ucI$ provided that
$\uP|_F=\ucI|_F$ for any face $F=\{x_i=1\}$ of $\Delta(r,n)$, $r>2$.

Assume, on the contrary, that a certain $P_\alpha\in\ucI$ is broken into pieces.
Choose a polytope $Q\subset\uP\cap P_\alpha$ such that the boundary of $Q$
contains the face $F=\{x_l=1\}\cap P_\alpha$, $l\not\in I_\alpha$.
A polytope $Q$ is realizable. In the corresponding configuration $D$,
the hyperplane $L_l$ is multiple (of multiplicity~$|I_\alpha^c|$).
and intersections of hyperplanes $L_j$, $j\in I_\alpha$ with $L_l$ are in general position
(because $F\subset Q$).
It follows that $D$ is central (except that $L_l$ is multiple).
If $Q\ne P_\alpha$ then there is at least one degeneracy, hyperplanes $L_j$, 
$j\in J\subset I_\alpha$, $|J|=r$
pass through a point $\beta\not\in L_l$.
Since not all hyperplanes $L_i$, $i\in I_\alpha$ pass through $\beta$,
there exist indices $k,k'\in J$ and $i\in I_\alpha$
such that a line $\cap_{i\in J\setminus\{k,k'\}}L_i$ intersects
$L_k$ and $L_{k'}$ at $\beta$, $L_l$ at $L_i$ at two other distinct points.
It follows that $\Delta_{\{k,k',l,i\},J\setminus\{k,k'\}}(2,4)$ is broken by 
$\uP$.
\end{proof}

\begin{proof}[Proof of~\eqref{KJFKJFKKHGD}.]
Let $M_C$ be as in~\eqref{KHGKHGJ} for a fixed lax hyperplane
arrangement $C$. Let $Z \subset \oo_{\ucI}$ be the fibre
over the point of $X_C(r,n)$ given by $C$, in the product
decomposition \eqref{sean3.6}. 

We consider lines $x_M:\,\bA^1 \to \cM(r,n)$, 
$x_M(z)=M_C+zM$ for $M\in \cM(r,n)$, and the induced regular map 
(which we abusively 
denote by the same symbol) $x_M: \, \bA^1 \to \oX_L(r,n)$.
We consider the limit of $x_M$ as $z \to 0$. 

We assume that $N$ has the lax order of~\eqref{sean3.8},
so for any $l$, points on $L_l$ of multiplicity greater than $r$ with respect to $L_1,\ldots,L_{l-1}$
are linearly independent. Let $p_l$ be the number of these points. 
A moments thought, and \eqref{sean3.6} yields the equality 
$$\sum_{\alpha \in \cI}(|I_a| - r)=\sum_i p_i = \dim(Z).$$

We now construct the columns of $M$. Suppose the first 
$l-1$ columns of $M$ are already constructed,
and consider column $l$. 
Let $e_i\in L_l$, $i=1,\ldots,p_l$ be points of multiplicity greater than 
$r$.
We include these $e_i$'s in the basis $e_1,\ldots,e_r$ and
write the $l$-th column in the dual basis. 
Let $V_i,W_i\subset N$, $i=1,\ldots,p_n$ be any choice of 
subsets as in~\eqref{KGLJGLKYF}
such that $V_i=\{i_1,i_2,i_3,l\}$, 
$|V_i\cap I^c_{e_{i}}|=1$, $W_i\subset I_{e_i}$. 

\begin{Claim}
For $i=1,\ldots,p_l$,
$$\lim\limits_{z\to0}CR_{V_i,W_i}(x_M(z))$$
does not depend on $M_{jl}$ for $j\ne i$
and depends non-trivially on $M_{il}$ (i.e. we can make this
limit any general value by varying $M_{il}$).
\end{Claim}
We can assume without loss of generality that $i_1\not\in I_{e_i}$
(otherwise take an appropriate automorphism of a cross-ratio function). 
Note the claim implies the result: First its clear that any
single choice of subsets $W,V$ as in \eqref{JFJFFJFJF} can be
chosen as $W_i,V_i$ for some $i$ and $l$. So all these cross
ratios are generic (i.e. take on values other than 
$\{0,1,\infty\}$) for general $M$. Now by \eqref{JFJFFJFJF} the
limit point is in $Z$. Now by the claim, we can vary
$\dim(Z)$ of the cross ratios completely independently by varying
$M$. Since $Z$ is smooth and connected by \eqref{sean3.6} it thus
follows that $Z \subset \oX_L(r,n)$ (set theoretically) 
and so since $C$ was arbitrary this completes the proof.

Let $W:=W_i$.
Then
$$\lim_{z\to0}CR_{V_i,W_i}(x_M(z))=
\lim_{z\to0}{\Det_{i_1i_2W}\Det_{i_3lW}\over\Det_{i_1i_3W}\Det_{i_2lW}}
$$
Notice that $\lim_{z\to0}\Det_{i_1i_2W}$ and $\Det_{i_1i_3W}$ are not 
zero - by assumption $L_{i_1}$
does not pass through $e_i$, but projections of any $r-1$ hyperplanes in 
$I_{e_i}$ from~$e_i$ are linearly independent.

So we have to demonstrate that
$$
\lim_{z\to0}{\Det_{i_3lW}\over\Det_{i_2lW}}
$$
does not depend on $M_{jl}$ for $j\ne i$
and depends not trivially on $M_{il}$.
Indeed, the constant terms of
$\Det_{i_3lW}$ and $\Det_{i_2lW}$ vanish, let's find coefficients at~$z$.
The $i$-th rows of the corresponding submatrices of $M_C$ is trivial, so
we can expand both determinants along this row and get
$$
\lim_{z\to0}{\Det_{i_3lW}\over\Det_{i_2lW}}=
{M_{ii_3}R_{ii_3}+M_{il}R_{il}+\ldots\over
M_{ii_2}R_{ii_2}+M_{il}R_{il}+\ldots},$$
where $R_{ij}$ are cofactors of the corresponding submatrices of $M_C$.
These cofactors are not trivial because projections of any $r-1$ 
hyperplanes in $I_{e_i}$ from~$e_i$ are linearly independent. So we see 
that the limit indeed does not depend on $M_{jl}$ for $j\ne i$ and is a 
M\"obius function
in $M_{il}$. This function can be made nontrivial by adding an  open condition
$M_{ii_3}R_{ii_3}\ne M_{ii_2}R_{ii_2}$.
\end{proof}

\begin{Proposition}\label{sean3.11}
Let $C$ be a lax configuration with multiple points $\cI$. If $|\cI|\ge2$ then
$\codim_{\oX_L(r,n)}\oo_{\ucI}\ge2$.
\end{Proposition}

\begin{proof}
We proceed by induction on $\sum|I_\alpha|$ using~\eqref{sean3.6}
and the following observation - a configuration
near the GIT-stable configuration is GIT-stable.
We will compare quantities
$\dim X_C(r,n)+\sum_{\alpha \in \cI}(|I_a| - r)$
for various configurations. Since all of them will be GIT stable,
we can substitute $X_C(r,n)$ by its $\PGL_r$-torsor $\bP_C(r,n)$,
the space of all configurations with prescribed multiplicities.

Assume first that there are some points of multiplicity greater than $r$.
Take the last hyperplane $L$ in the lax order that contains such point.
Move $L$ a little bit to take it off this point but keep all other points 
of multiplicity greater
than $r$ on $L$ (this is possible because they are linearly independent).
If there are no such points, keep some point of multiplicity $r$ 
(if there are any of them on $L$).
Then the dimension of the configuration space will increase by at least one,
the sum $\sum |I_\alpha|-r$ will decrease by at most one and $|\cI|$ is 
still at least $2$.
At the end, there will be at least two points $A$, $B$ of multiplicity 
$r$ and no points of higher multiplicity.
Now take a hyperplane through $A$ and move it
keeping $B$ if it belongs to this hyperplane.
This will increase the dimension of the configuration space by at least one
but the result will still not be generic, thus having codimension at least one.
\end{proof}

\begin{proof}[Proof of \eqref{:LKJH:KJH}] From \eqref{sean3.6} the
generic stratum of $B_{I_a}$ is smooth and connected, and
codimension one in $U$. By \eqref{sean3.11} all other boundary strata of
$U$ are lower dimensional. The boundary of $U$ is (by definition)
the scheme-theoretic inverse
image of the boundary of $U_{\ucI}$, and so Cartier, and 
in particular pure codimension one, by 
\eqref{sean3.4}.
It follows the $B_{I_a}$ are
irreducible, Cartier, and their union is the full
boundary. They are generically smooth by \eqref{sean3.6}. 
The proof of \eqref{sean1.8} shows that their complement, the
main stratum $X(r,n)$ is isomorphic to an open subset
of affine space, and thus has trivial divisor class group.
Thus $U$ is factorial. Now it is smooth generically along
the Cartier divisors $B_{I_a}$ by \eqref{sean3.6}.
In the open set $U_{\ucI} \subset \cA$ the stratum
$\cA_{\ucI}$ is the scheme-theoretic intersection of the
boundary divisors that contain it (this is true in any
toric variety). Thus $\oo_{\ucI}$ is scheme-theoretically
the intersection of
the boundary divisors of $U$ that contain it.
\end{proof}

Proof of~\eqref{LKHGKKG} now follows from \eqref{sean3.6} and 
\eqref{sean3.13} below.

\begin{Proposition}\label{sean3.13}
Let $X$ be a normal variety.
Let $B_i$ be irreducible $Q$-Cartier Weil divisors.
If $K_X + \sum B_i$ is log canonical, then the intersection
$$
B_1 \cap B_2 \dots \cap B_n
$$
is (either empty or) pure codimension $n$.
\end{Proposition}

\begin{proof} We can intersect with a general hyperplane
to reduce to the case when $n$ is the dimension of $X$
and then apply \cite[18.22]{FA}.
\end{proof}

\section{The membrane}\label{sean4}

\begin{Review}\label{sean4.1}\textsc{Tits buildings.}
We begin with some background on buildings. For further details
see \cite[\S 1]{Mustafin78}.
$R = k[[z]]$ is the ring of power series and $K$ is its quotient field. $k$,
as throughout the paper, is an algebraically closed field.
Let $V = k^r$ and $V_K = V \otimes_k K$.

We follow \cite[pg 108]{Spanier} for elementary definitions and
properties of simplicial complexes. In particular for us
a simplicial
complex is a set (the vertex set of the complex)
together with a collection
of finite subsets called simplicies.
\end{Review}

\begin{Definition}
The {\it total Grassmannian} $\Gr(V)$ is a simplicial complex
of dimension $r-1$. Its vertices are non-trivial subspaces of $V$
(in particular $V$ itself is a vertex). A collection of distinct subspaces
forms a simplex iff they are pairwise incident (i.e. one is contained
in the other), from which it easily follows that $m-1$ simplicies
correspond to flags of non-trivial subspaces
\begin{equation}\label{simplexGr}
0=U_0 \subset U_1 \subset U_2 \dots \subset U_m.\cooltag
\end{equation}

The compactified {\em affine Tits building} $\ocB$ is a simplicial
complex of dimension $r-1$, with vertices given by equivalence classes
of non-trivial free $R$-submodules of $V_K$, where $M_1,M_2$ are equivalent
iff $M_1 = cM_2$ for some $c \in K^*$.
A collection of distinct equivalent classes is a simplex iff they are pairwise incident,
where $[M_1]$ and $[M_2]$ are called {\it incident} if after rescaling
$$
zM_1 \subset M_2 \subset M_1.
$$
Incidence is easily seen to be symmetric.
\end{Definition}

\begin{Lemma}\label{sean4.1.1}
Distinct pairwise incident classes
$\Gamma_1,\dots,\Gamma_m \in \ocB$ form a $m-1$ simplex $\sigma$ iff
after reordering
there exist representatives $[M_i] = \Gamma_i$ so that
\begin{equation}\label{simplexsigma}
zM_m = M_0 \subset M_1 \subset M_2 \dots \subset M_m.\cooltag
\end{equation}
\end{Lemma}

\begin{proof} \cite[1.1]{Mustafin78}.
\end{proof}

Note that we can rescale so that any $M_i$ in the simplex
is in the position of $M_m$.

The affine Tits building $\cB\subset \ocB$ is the
full subcomplex of equivalence
classes of lattices, i.e. free submodules of rank $r$
(full subcomplex means
that a subset of $\cB$ forms a simplex iff it
does in $\ocB$). As a set, $\cB=\PGL_r(K)/\PGL_r(R^*)$.

\begin{Review}\textsc{Stars.}
Recall that the {\it Star} of a simplex $\sigma$ in a simplicial
complex $C$ is a subcomplex
$$
\Star_{\sigma} C = \bigcup_{\sigma \subset \partial \sigma'} \sigma'.
$$
Notice that if $\tau\subset\partial\sigma$ then  $\Star_\sigma$ is canonically
a subcomplex of $\Star_\tau$.
\end{Review}

\begin{Lemma}\label{LKJGLKHG}
For any lattice $[\Lambda]\in\cB$, there exists a canonical isomorphism
$$\Star_\Lambda\cB\simeq\Gr(\oLambda),$$
where $\oLambda=\Lambda/z\Lambda$.
More generally, for any simplex $\sigma$~\eqref{simplexsigma} in $\cB$,
$$\Star_\sigma\cB\simeq\Gr(M_m/M_{m-1})*\ldots*\Gr(M_1/M_{0}),$$
where $*$ denotes the join of simplicial complexes.
\end{Lemma}

\begin{proof} Indeed, lattices between
$z\Lambda$ and $\Lambda$ are obviously in incidence preserving
bijection with subspaces of $\oLambda$. Similarly,
lattices that fit in the flag~$\sigma$ correspond bijectively
to subspaces in one of $M_i/M_{i-1}$.
\end{proof}

It is probably worth mentioning that for any simplex $\sigma$ of $\Gr(V)$ as in \eqref{simplexGr},
$\Star_\sigma\Gr(V)=\Gr(U_m/U_{m-1})*\ldots*\Gr(U_1/U_0)$.

\begin{Review}\label{sean4.2}\textsc{Retractions.}
For a lattice $\Lambda$, and a non-trivial subset $\Theta\subset V_K$
(e.g.~an $R$-submodule or an element)
we can find unique $a > 0$ so that $z^a \Theta \subset \Lambda$,
$z^a \Theta\not \subset z\Lambda$. We define $\Theta^\Lambda := z^a \Theta \subset \Lambda$.
We denote
$$
R_{\Lambda}: \ocB \rightarrow \Star_\Lambda\cB
$$
the map that sends a submodule $M$ to $M^{\Lambda}+z\Lambda$.

\begin{Lemma}\label{sean4.2.1}
$R_{\Lambda}$ is a retraction map of simplicial complexes.
\end{Lemma}

\begin{proof} It's clear the map preserves incidence
and is identity on $\Star_\Lambda\cB$.
\end{proof}

More generally,
for any simplex \eqref{simplexsigma} in $\cB$
there is a retraction
$$
R_{\sigma}: \ocB \rightarrow \Star_\sigma\cB
$$
defined as follows: $R_\sigma(M)=M^{\Lambda_k}+\Lambda_i$,
where $i$ is maximal such that $M^{\Lambda_k}\not\subset\Lambda_i$.
We denote by $\Res_\Lambda$ and $\Res_\sigma$ compositions of $R_\Lambda$
and $R_\sigma$ with isomorphisms of \eqref{LKJGLKHG}.
\end{Review}

\begin{Review}\label{sean4.3}\textsc{Convex structure.}
A subset of $\ocB$ is called {\it convex} if it is closed under
finite $R$-module sums. For a collection of free $R$-submodules
$\{M_{\alpha}\}$
their {\it convex hull} in $\ocB$, denoted $[M_{\alpha}]$
is the subcomplex with vertices (with representatives) of form
$\sum c_{\alpha} M_{\alpha}$, $c_{\alpha} \in K$.
This is obviously
the smallest convex subset that contains all the $[M_{\alpha}]$.

A subset of $\Gr(V)$ is called convex if it is closed under
finite sums (of subspaces). We write $[U_{\alpha}]$ for the convex hull
of subspaces $\{U_{\alpha}\}\subset \Gr(V)$.

It is clear that stars of simplices in $\cB$ and $\Gr(V)$ are convex and that
isomorphisms~\eqref{LKJGLKHG} preserve the convex structure.
\end{Review}

\begin{Example}
For $T=\{g_1,\dots,g_r\}$ a basis of $V_K$, the convex hull
$[T]$ is called an {\it apartment}.
It is the
set of equivalence classes $[M]$ such that $M$ has an
$R$-basis $c_1g_1,\dots,c_r g_r$ for some $c_i \in K$.
\end{Example}

Retractions commute with convex hulls:

\begin{Proposition}\label{sean4.3.1}
Let $\{M_{\alpha}\} \subset \ocB$ be
a subset and $\sigma \subset \cB$ a simplex. Then
$$
[R_{\sigma} M_{\alpha}] = R_{\sigma}[M_{\alpha}].
$$
If $\sigma \subset [M_{\alpha}]$ then both sides are also equal to
$\Star_{\sigma} [M_{\alpha}]$.
\end{Proposition}

\begin{proof} We leave it as an exercise to the reader.
\end{proof}

\begin{Lemma}[{\cite{Faltings01}}]\label{sean4.3.2}
The convex hull of a finite subset of $\cB$ is finite.
\end{Lemma}

\begin{Review}\label{sean4.4}\textsc{Membrane.}
Let
$$
\cF := \{Rf_1,\dots,Rf_n\}
$$
be a collection of $n$ rank $1$ submodules, such that
$[f_1,\ldots,f_n]=V_K$.
The convex hull $[\cF]\subset\ocB$ we call the {\it membrane}.
\end{Review}

\begin{Lemma}\label{sean4.4.1}
$[\cF]$ is the union of apartments
$[T]$ for subsets $T \subset \cF$, $|T| = r$.
\end{Lemma}

\begin{proof} This is immediate from the definitions, and Nakayama's
lemma.
\end{proof}

\begin{Review}\label{sean4.5}\textsc{Membranes as tropical subspaces.}
We begin by recalling the construction of the tropical variety
(also called a non-Archimedean amoeba or a Bieri--Groves set),
see \cite{SS} for details.
Let $H=\bG_m^n$ be an algebraic torus.
Let $\overline K$ be the field of generalized Puiseaux series
$\sum\limits_{\alpha\in I\subset\bR} c_\alpha z^\alpha$,
where $I$ is a locally finite subset of $\bR$, bounded below
(and is allowed to vary with the series).
There is an evaluation map
\def\ord{{\text{\rm ord}}}\def\Ord{{\text{\rm Ord}}}
$$\ord:\,H(\overline K)\to H(\overline K)/ H(\oR)=\bR^n$$
where $\oR \subset \oK$ is the subring of series for which
$I \subset \bR_{\ge0}$.
For any subvariety $Z\subset H$, $\ord(Z)$
is called the {\em tropicalisation} of $Z$.
It is a polyhedral complex of dimension $\dim Z$.
If $Z$ is invariant under dilations
then $\ord(Z)$ is invariant under diagonal translations and
we consider
$$\Ord(Z)=\ord(Z)\mod \bR(1,\ldots,1).$$

Let $\cF := \{Rf_1,\dots,Rf_n\}$ be as in \eqref{sean4.4}.
Consider the map
$$\Phi:\,V_K^\vee\to K^n,\quad F\mapsto (F(f_1),\ldots,F(f_n)),$$
Let $Z=\Phi(V_K^\vee)\cap H$. Then $Z$ is of course the intersection
with $H$ of the $r$-plane
spanned by the rows of the $r \times n$ matrix with columns
given by $f_i$'s.
Its tropicalisation
$\Ord(Z)\subset\bR^{n-1}$ is called {\em a tropical projective subspace}.

For any simplicial complex $C$, we denote by $|C|$ the corresponding
topological space  (obtained by gluing physical simplicies).
Recall that $|\cB|$ can be identified with the space of equivalence classes
of additive norms on $V_K$, where an additive norm $N$ is a map $V_K(K)\to\bR\cup\{\infty\}$
such that
$$N(cv)=\ord(c)+N(v)\quad \text{for any}\quad c\in K,\ v\in V_K(K),$$
$$N(u+v)\ge\min(N(u),N(v))\quad \text{for any}\quad u,v\in V_K(K),$$
and
$$N(u)=\infty\quad \text{iff}\quad u=0.$$
Two additive norms are equivalent if they differ by a constant.
For a norm $N$ let
$\tPsi(N) = (N(f_1),\ldots,N(f_n)) \in \bR^n$.
Now consider
$$
\Psi:\,|\cB|\to\bR^{n-1},\quad \Psi([N])= \tPsi(N) \mod
\bR(1,\ldots,1).
$$
The map is continuous because the topology on $|\cB|$ is exactly the
topology
of point-wise convergence of norms.
The following theorem is our version of the tropical Gelfand--Macpherson transform.
\end{Review}

\begin{Theorem}\label{sean4.5.1}
$\Psi$ induces a homeomorphism
$|[\cF]|\simeq\Ord(Z)$.
\end{Theorem}

\begin{proof}
Let $\Omega$ be the
Drinfeld's symmetric domain -- the complement
to the union of all $K$-rational hyperplanes in $V^\vee_K(\overline K)$.
There is a surjection \cite{Drinfeld}
$$D:\,\Omega\to|\cB|,\quad F\mapsto[v\to\ord\, F(v)]\quad \text{for}\
v\in V_K(K),$$
here we interpret $|\cB|$ as the set of equivalence classes of norms.
The following diagram is obviously commutative:
$$\begin{matrix}
\Omega&\mathop{\longrightarrow}\limits^D&|\cB|\cr
\phantom{\!\!\scriptstyle\Phi}\downarrow{\!\!\scriptstyle\Phi}&&
\phantom{\!\!\scriptstyle\Psi}\downarrow{\!\!\scriptstyle\Psi}\cr
H(\overline K)&\mathop{\longrightarrow}\limits^{\Ord}&\bR^{n-1}
\end{matrix}$$
\def\Im{\text{Im}}
It follows that $\Im(\Psi)\subset\Ord(Z)$.

For any lattice $\Lambda\in\cB$, the corresponding norm $N_\Lambda$
is as follows:
$$N_\Lambda(v)=\{-a\,|\,z^a v\in\Lambda\setminus z\Lambda\}\in\bZ.$$
In particular, $\Psi(\cB)\subset\bZ^{n-1}$.
Also, it follows easily from definitions that $\Psi$ is affine on simplicies of $|\cB|$
and uniformly continuous.
Since $\Ord(Z)$ is a polyhedral complex, it remains to check that
for any $\bQ$-point of $\Ord(Z)$, there exists a unique
$\bQ$-point of $|[\cF]|$ that maps onto it
(a $\bQ$-point of $|\cB|$ means a point of some simplex with rational barycentric
coordinates). Now we can pass from $K$ to Puiseaux series $k[[z^{1/m}]]$
with sufficiently
large $m$ (this does not change neither $\Ord(Z)$ nor $|[\cF]|$, see
also \eqref{sean6.3} for another version of this baricentric trick) and
it remains to check the latter statement for $\bZ$-points.
Substituting $f_i$'s by $z^{a_i}f_i$'s, we can assume that this point
is $O=(0,\ldots,0)$. Now we claim that if $O\in\ord(Z)$ then
$\tPsi(\Lambda)=O$ for $[\Lambda] \in[\cF]$
if and only if $\Lambda=Rf_1+\ldots+Rf_n$.

Suppose $\tPsi(Rf_1+\ldots+Rf_n)\ne O$, i.e.~$f_j\in z (Rf_1+\ldots+Rf_n)$ for some~$j$.
By Nakayama's lemma, we can assume without loss of generality that $Rf_1+\ldots+Rf_n=Rf_1+\ldots+Rf_r$.
Therefore,
$f_j\in z (Rf_1+\ldots+Rf_r)$. But then for any $F\in V_K^\vee(\overline
K)$,
if $\ord\, F(f_i)=0$ for $i\le r$ then $\ord\, F(f_j)>0$. But this
contradicts
$O\in\ord(Z)$.

Now take any lattice
$\Lambda=Rz^{a_1}f_1+\ldots+Rz^{a_n}f_n$.
We can assume without loss of generality that
$\Lambda=Rz^{a_1}f_1+\ldots+Rz^{a_r}f_r$. Let $\tPsi(\Lambda)=O$.
Then $f_i\in\Lambda$ for any~$i$,
therefore, $\Lambda=Rf_1+\ldots+Rf_n$.
\end{proof}

\begin{Review}\label{sean4.5.2}\textsc{A generalization of the visible contour family.}
Let
$$H=\bG_m^{n-1}(K) \subset \bP^{n-1}(K).$$
Start with any subvariety $Z \subset H$.
Consider the point
$$
[\oZ] \in\Hilb_{\bP^{n-1}}(K).
$$
Assume for simplicity that no element of
$H$ preserves $Z$.
We consider the orbit closure
$\overline{H \cdot [\oZ]} \subset\Hilb_{\bP^{n-1}}(K) $, a toric variety for
$H$, which gives us
a $K$-point
$$\left[\overline{H \cdot [\oZ]}\right] \in \Hilb(\Hilb_{\bP^{n-1}}),$$
which we
can think of as a one parameter family of toric varieties over $k$.
We can then let $z$ approach zero and degenerate the toric varieties
to a {\it broken} toric variety in $\Hilb_{\bP^{n-1}}$.
Note if we
start with $\oZ$ a general $r-1$ plane, then
$Z \subset \oZ$ is the complement of $n$ general hyperplanes
in $\bP^{r-1}$, the corresponding component of $\Hilb_{\bP^{n-1}}$ is
$G(r,n)$, and the degeneration takes place in Kapranov's
family of {\it broken} toric varieties,
$\cT \rightarrow G(r,n)//H$ of \eqref{LIUYFKJYF}.
The visible contour family also
generalizes: Since $Z \subset H$ we can embed $Z$ in
$H \cdot [\oZ]$ by
$t \rightarrow t^{-1} \cdot [\oZ]$, so the image has the strange
expression $Z^{-1} \cdot [\oZ]$. We can then consider
$$
\left[\overline{Z^{-1} \cdot [\oZ]}\right] \in \Hilb(\Hilb_{\bP^{n-1}})
$$
Note $Z^{-1} \cdot [\oZ] \subset \Hilb_{\bP^{n-1}}$ is precisely
the set of points which (thought of as subschemes of $\bP^{n-1}$)
contain the point $(1,\dots,1)$, so in the case of a general
linear space $\oZ$ this is precisely Kapranov's visible contour, \eqref{sean2.6}.

In the linear case the broken toric variety is described by
a matroid decomposition of $\Delta(r,n)$ which reflects the
combinatorics of the simplicial complex $[\cF]$ and so by
\eqref{sean4.5.1}, of the tropical variety $\Ord(\oZ)$. It is natural
to wonder if this holds in general.
\end{Review}

\input part2.tex

%% file: part2.tex
\section{Deligne Schemes}\label{sean5}For convex
$$
\Stab \subset Y \subset \cF
$$
and $p:\bS_Y \rightarrow \Spec(R)$ the Deligne scheme, we
have by \eqref{sean5.14} the natural vector bundle $\Omega^1_p(\log \bB)$,
see \ref{sean9}.

\begin{Review}
Now we turn to the proofs of \eqref{sean1.4} and \eqref{sean1.16}--\eqref{sean1.18}.
We follow the notation of the introduction and \ref{sean4}.
Here we prove the pair $(\bS_Y,S_Y + \bB)$ of \eqref{sean1.4} has normal
crossings, \eqref{sean5.14}. Global generation is considered in \ref{sean6}.
\end{Review}

\begin{Review}\textsc{Deligne Functor \cite{Faltings01}.}
Let $Y \subset \cB$ be a finite set.
A {\em Deligne functor}  $\bS_Y$ is a
functor from $R$-schemes to sets,
a $T$-valued point $q$ of which
consists of a collection of equivalence
classes of line bundle quotients
$$
q_M:\, M_T \twoheadrightarrow L(M_T)
$$
for each lattice $[M] \in Y$, where $M_T := T \times_{R} M$,
where two quotients are equivalent if they have the same
kernel, satisfying the compatibility requirements:
\begin{itemize}
\item
For each inclusion $i:\,N \hookrightarrow M$,
there is a commutative diagram
$$
\begin{CD}
N_T  @>{q_N}>> L(N_T)\\
@Vi_TVV            @VVV \\
M_T @>{q_M}>> L(M_T)
\end{CD}
$$
\item
Multiplication by $c\in K^*$ gives an isomorphism
$$
\begin{CD}
\ker q_M @>\cdot c>> \ker q_{cM}
\end{CD}
$$
\end{itemize}
It is clear from this definition that $\bS_Y$ is represented by a closed subscheme
$$\bS_Y\subset\prod\limits_{[M]\in Y}\cP(M),$$
$(\bS_Y)_K=\cP(V_K)$, and $\bS_Y$ contains the Mustafin's join~\eqref{LKUJHLKGKHG}.
\end{Review}

\begin{ETheorem}\textsc{Theorem \cite{Faltings01}.}\label{sean1.14}
Assume $Y$ is finite and convex. Then
$\bS_Y$ is smooth and irreducible (in particular it is
isomorphic to the  Mustafin's join).
Its special fiber $S_Y=(\bS_Y)_k$ has normal crossings.
\end{ETheorem}

We begin by explaining Faltings' proof of \eqref{sean1.14}, recalling
and expanding upon the three paragraphs of
\cite[pg. 167]{Faltings01}. This is the substance
of \eqref{sean5.2}--\eqref{sean5.9.1}.
For this $Y \subset \cB$ is an arbitrary
finite convex subset.
Beginning
with \eqref{sean5.10} our treatment diverges from \cite{Faltings01}.
We specialize to convex subsets $Y \subset [\cF]$
as in \eqref{sean1.4} and consider singularities of the natural boundary.

\begin{Review}\label{sean5.2}\textsc{Maximal Lattices \cite{Faltings01}.}
Consider a $k$-point of $\bS_Y$, i.e. a compatible
family of one dimensional $k$-vector space quotients
$$
q_M: \oM\twoheadrightarrow L(\oM),\quad [M] \in Y,
$$
where
$$\oM=M/zM=M\otimes_Rk.$$
This gives a partial order on $Y$: $[N]\le_q[M]$ iff the composition
$$
\oN=\overline{N^M}\to\oM\twoheadrightarrow L(\oM)
$$
induced by inclusion $N^M\subset M$
is surjective, and thus by compatibility, canonically
identified with $q_N$.
In this case we also say that $q_M$ does not vanish on $N$.

A lattice $[M] \in Y$ is called {\it maximal} for $q$
if it is maximal with respect to the order $\le_q$. In other words,
$q_N$ vanishes on $M$ for any $[N]\in Y$, $[N]\ne[M]$.
Since $Y$ is finite, it follows that
for each $[N] \in Y$ there exists a maximal lattice
$[M] \in Y$ such that $[N]\le_q[M]$.
\end{Review}

\begin{ETheorem}\textsc{Lemma \cite{Faltings01}.}\label{sean5.2.2}
Maximal lattices are pairwise incident.
\end{ETheorem}

It follows by \eqref{sean4.1.1} that maximal lattices form a simplex~$\sigma$~\eqref{simplexsigma}.

\begin{Corollary}\label{sean5.2.3}
A  $k$-point of $\bS_Y$
with a simplex of maximal lattices $\sigma$ is equivalent to a collection of hyperplanes
$$
H_i \subset M_{i}/M_{i-1},\quad m \geq i \geq 1,
$$
which do not contain $\Res_{\sigma}[M]$ for any $[M] \in Y$.
\end{Corollary}

\begin{proof}
Choose a $k$-point $q$ of $\bS_Y$
with a simplex of maximal lattices $\sigma$.
Let $[M]\in Y$ and let $i$ be such that
$\Res_\sigma[M]\in\Gr(M_i/M_{i-1})$.
Rotating $M_i$'s if necessary,
we can assume that $i=1$.
Let $M_i$ be a maximal lattice such that $[M]\le_q[M_i]$.
Then $q_{M_i}$ does not vanish on~$M$, and therefore does not
vanish on $M_1$. But $M_1$ is maximal, so $i=1$.
It follows that hyperplanes $H_i=(\ker q_{M_i})/M_{i-1}$
don't contain $\Res_{\sigma}[M]$ for any $[M] \in Y$.
And it is clear that these hyperplanes determine $q$.
\end{proof}

\begin{Definition}\label{sean5.3}
Let $[M] \in Y$. We let
$\tqp(\oM) \subset \bS_Y$ be the subfunctor of
compatible quotients such that
for each $[N] \in Y$, $N=N^M$,
the quotient
$q_{N}:N_T \twoheadrightarrow L(N_T)$
vanishes on $(N\cap zM)_T\subset N_T$.
\end{Definition}

It's clear $\tqp(\oM)$ is represented by a closed subscheme
of $\bS_Y$.

\begin{Lemma}\label{sean5.4}
The $k$-points of
$\tqp(\oM) \subset \bS_Y$ are precisely the set of $k$-points
of $\bS_Y$ for which $M$ is a maximal lattice.
\end{Lemma}

\begin{proof} Consider a $k$-point $q$ of $\tqP(\oM)$.
Suppose $M \subsetneq N$,
$[N] \in Y$. Then $N^M = z^{k} N$, for
some $k >0$.
$q_{z^kN}$ vanishes on $z^k N \cap zM$,
by the  definition of $\tqp(\oM)$.
So by compatibility of quotients with scaling
$q_N$ vanishes on $N \cap z^{1-k}M$, which contains $M$.
Thus $M$ is a maximal lattice.

Conversely, suppose $M$ is
maximal for a $k$-point $q$. Take $[N] \in Y$ such that $N=N^M$.
By maximality $q_{z^{-1}N + M}$ vanishes
on $M$, thus by compatibility, $q_{N + zM}$ vanishes
on $zM$, thus, again by compatibility,  $q_N$ vanishes on  $N\cap zM$. So the point lies in
$\tqP(\oM)$.
\end{proof}

\begin{Definition}
Let $V$ be a finite-dimensional $k$-vector space,
and let $\cW \subset \Gr(V)$ be a finite convex collection
of subspaces that includes $V$.  Let
$\bl(\qP(V),\cW)$ be the functor from
$k$-schemes to sets which assigns to each $T$ the
collection of line bundle quotients
$W_T \twoheadrightarrow L(W_T)$, $W \in \cW$, $W_T$ the
pullback,  compatible with the inclusion maps between the $W$,
i.e. the composition
$$\begin{CD}
A_T \rightarrow B_T @>{q_B}>> L(B_T)\end{CD}
$$
factors through $q_A: A_T \twoheadrightarrow L(A_T)$
for $A \subset B$, $A,B \in \cW$.
\end{Definition}

\begin{Proposition}\label{sean5.6}
There is a canonical identification
$$
\tqP(\oM) = \bl(\qP(\oM),\Res_{M}(Y)).
$$
\end{Proposition}

\begin{proof}
Immediate from the definitions.
\end{proof}

\begin{Proposition}\label{sean5.5}
$\bl(\qP(V),\cW)$ is represented by the
closure of the graph of the product of canonical rational maps
$\qP(V) \dasharrow \qP(W)$, $W \in \cW$. Furthermore
$\bl(\qP(V),\cW)$ is smooth.
\end{Proposition}

\begin{proof} We induct
on the number of subspaces in $\cW$. When $\cW = \{V\}$
the result is obvious.
In any case it is
clear the functor is represented
by a certain closed subscheme
$$
X \subset \prod_{W \in \cW} \qP(W).
$$
Let  $\qP^0(V)\subset\qP(V)$ be an open subset of quotients that don't vanish
on any $W\subset\cW$.
Then
$\qP^0(V)$ is an open subset of $X$, its closure $X'$ in $X$
is the closure of the graph
in the statement.

Take a fixed closed point $q^0$. We will show that $X=X'$ near $q^0$.
Let $W \in \cW$  be a maximal subspace such that $q^0_V$ vanishes on $W$.
If there are none then $q^0\in\qP^0(V)$ and so $X=X'$ near $q^0$.

Let $D \subset X$ be the subscheme of compatible quotients
so that $q_V$ vanishes on $W$, let $D^0\subset D$ be a sufficiently small neighbourhood of $q^0$.
Let $\cW_W \subset \cW$
be those subspaces contained in $W$, clearly $\cW_W$ is
convex.
Take any $E \in \cW$ and $q\in D^0$.
If $E \not\in\cW_W$, then $q^0_V$ (and hence $q_V$) does not vanish on~$E$,
from which it follows that
$E \rightarrow L(V)$ is surjective, and thus identified with
$q_E$. It follows easily that
$D^0$ is represented by an open subset of $\bl(\qP(W),\cW_W) \times \qP(V/W)$. In particular,
by induction, $D^0$ is connected and smooth of dimension
$\dim V -2$.

Claim: $D^0 \subset X'$. As $D^0$ is integral it's enough to check this
on some open subset of $D^0$. We consider the open subset
where $q_W$ does not vanish on any $E \in \cW_W$, and
$q_V$ does not vanish on any $E \not \in \cW_W$. This is
naturally identified with an open subset for
$\cW =\{V, W\}$, and so we reduce to this case. But in this case
it is easy to see that
$X = X'$ is the blowup of $\qP(V)$ along $\qP(V/W)$
and so obviously $D^0\subset X'$.

By the Claim $X$ has dimension
at least $\dim V-1$ along $D^0$. $D \subset X$ is locally
principal, defined by the vanishing of a map between the
universal quotient line bundles for $W$ and $V$. It
follows from \eqref{sean5.5.1} below that $X$ is smooth, and
equal to $X'$ along $D^0$.
\end{proof}

The following is well known:

\begin{Lemma}\label{sean5.5.1}
Let $(A,m)$ be a local Noetherian
ring of Krull dimension at least $d$. Assume
$A/f$ is regular of dimension at most $d-1$ for $f \in m$. Then
$A$ is regular of dimension $d$. \end{Lemma}
\begin{proof}
$$
\dim m/m^2 \leq 1 + \dim m/(m^2 + f) = 1 + \dim A/f = d
\leq \dim A \leq \dim m/m^2.
$$
\end{proof}

\begin{Definition}
Define the {\it depth} of $W \in \cW$ to be the
largest $d \geq 0$ so that there is a proper flag
$$
W = W_0 \subset W_1 \dots \subset W_d = V
$$
with $W_i \in \cW$. Let $\cW_{\le m} \subset \cW$ be
the subset of subspaces of depth at most~$m$.
Let $\cW_m \subset \cW_{\le m}$ be the
subset of subspaces of depth exactly~$m$.
\end{Definition}

Notice that $\bl(\qP(V),\cW_{0})=\qP(V)$ and
$\bl(\qP(V),\cW_{\le N})=\bl(\qP(V),\cW)$  for $N\gg0$.
Thus the next proposition shows that the canonical map
$\bl(\qP(V),\cW)\to\qP(V)$ is an iterated blowup
along smooth centers.

\begin{Proposition}\label{sean5.5.2}
A forgetful functor
\begin{equation}\label{KJFKFFK}
p: \bl(\qP(V),\cW_{\le m+1}) \rightarrow \bl(\qP(V),\cW_{\le m})\cooltag
\end{equation}
is represented by the blowup along the union of the
strict transforms of $\qP(V/W)\subset\qP(V)$ for $W \in \cW_{m+1}$
(which are pairwise disjoint).
\end{Proposition}

\begin{proof}
Let $W \in \cW_{m+1}$. We claim that the
strict transform of $\qP(V/W)$ represents
the subfunctor $X_W$ of
$\bl(\qP(V),\cW_{\le m})$ of compatible quotients such
that $q_E$ vanishes on $E \cap W$ for all
$E \in \cW_{\le m}$. This subfunctor is naturally identified with
$\bl(\qP(V/W),\cW_{\le m}^W)$, where
$\cW_{\le m}^W$ is the (obviously convex) collection
of subspaces $(E + W)/W \subset V/W$, for
$E \in \cW_{\le m}$.
By \eqref{sean5.5}, it is smooth and connected, and thus
the strict transform. For disjointness: if $W',W''\in\cW_{m+1}$ then
$\tilde W:=W' + W'' \in \cW_{\le m}$ and it is not
possible for $q_{\tilde W}$ to vanish both on $W'$ and
on $W''$, thus the strict transforms are disjoint.

The map \eqref{KJFKFFK}  is obviously an isomorphism
outside the union of subfunctors $X_W$. Take
$W \in \cW_{m+1}$. The inverse image
$$
p^{-1}(X_W) \subset \bl(\qP(V),\cW_{\le m+1})
$$
is naturally identified with
$\qP(W) \times \bl(\qP(V/W),\cW_{\le m}^W)$. In particular
by \eqref{sean5.5} it is a smooth connected Cartier divisor. It follows that the
exceptional locus of $p$ is the disjoint
union of these divisors.
Thus $p$ factors through the proscribed blowup, and
the induced map to the blowup will have no exceptional
divisors and is thus an isomorphism (as domain and
image are smooth).
\end{proof}

\begin{Definition}
For a subset
$\sigma \subset Y$, consider the intersection
$$
\tqP(\osigma) := \bigcap_{M \in \sigma} \tqP(\oM) \subset \bS_Y.
$$
\end{Definition}

\begin{Proposition}\label{sean5.7}
$\tqP(\osigma)$
is non-empty iff $\sigma$ is a simplex~\eqref{simplexsigma}.
Consider the convex subset $\Res_{\sigma}(Y)$,
a collection of convex subsets $\cW_{i} \subset \Gr(M_i/M_{i-1})$. There is
a canonical identification
$$
\tqP(\osigma) = \prod_{m \geq i \geq 1}
 \bl(\qP(M_i/M_{i-1}),\cW_i) =: \bl(\qP(\osigma),
\Res_{\sigma}Y).
$$
\end{Proposition}

\begin{proof} By \eqref{sean5.4} the $k$-points of the intersection are
exactly those for which all $M \in \sigma$ are maximal. Thus
if it is non-empty, $\sigma$ is a simplex by \eqref{sean5.2}.
The expression for the intersection
is immediate from the definition of $\Res$ (see \eqref{sean4.2}),
and the functorial definitions of $\tqP(\oM)$ and $\bl$.
\end{proof}

\begin{Remark}\label{sean5.7.1}
Observe by \eqref{sean5.6}--\eqref{sean5.7} that the
special fibre $S_Y$ has normal crossings. Moreover it
can be canonically defined purely in terms of the subcomplex
$Y \subset \cB$. Indeed by \eqref{sean5.6} its irreducible components
and their intersections are encoded by the
$\bl(\qP(\osigma),\Res_{\sigma}(Y))$ for simplicies
$\sigma \subset Y$, and
by \eqref{sean4.3.1} we have canonical identifications
$$
\Res_{\sigma} Y = \Star_{\sigma} Y \subset
\Star_{\sigma} \cB = \Star_{\osigma} \Gr(\oLambda_m).
$$
\end{Remark}

\begin{Definition}
Let $\sigma \subset Y$ be a simplex. Let
$U(\sigma) \subset \bS_Y$
be the open subset whose complement is the closed subset
of the special fibre given by the union of $\tqp(\oN)$,
$[N] \in Y \setminus \sigma$.
\end{Definition}

\begin{Lemma}\label{sean5.8}
$U(\sigma)$ is the union of the generic fibre together with
the open subset of the special fibre consisting of all
$k$-points whose simplex of maximal lattices \eqref{sean5.2}
is contained in $\sigma$. It represents the following subfunctor:
Let $\sigma$ be the simplex \eqref{simplexsigma}. For $[M] \in Y$
choose minimal $i$ so that
$M^{M_m} \subset M_i$.
A $T$-point of $\bS_Y$ is a point of $U(\sigma)$ iff
the composition
$$
M_T \rightarrow (M_i)_T \twoheadrightarrow L((M_i)_T)
$$
is surjective for all $[M] \in Y$.
\end{Lemma}

\begin{proof} Immediate from \eqref{sean5.4} and the definitions.
\end{proof}

Note by \eqref{sean5.2} that the $U(\sigma)$ for $\sigma \subset Y$
give an open cover of $\bS_Y$.
Faltings proves $U(\sigma)$ is non-singular, and
semi-stable over $\Spec(R)$, by writing down explicit local
equations, \cite[pg 167]{Faltings01}.
This can also be seen from the following:

\begin{Proposition}\label{sean5.9}
Let $\sigma \subset Y$ be the simplex \eqref{simplexsigma}.
Let $U \subset \qP(M_m)$ be the open subset of quotients
$M_m \rightarrow L$ such
that $N^{M_m} \rightarrow L$ is surjective for all
$[N] \in Y \setminus \sigma$.

Let $q:\bl(\qP(M_m),\osigma) \rightarrow \qP(M_m)$ be the iterated blowup
of $\qP(M_m)$ along the flag of subspaces of its
special fibre
$$
\qP(M_{m}/M_{m-1}) \subset \qP(M_{m}/M_{m-2}) \dots
\subset \qP(M_m/M_1) \subset \qP(M_m/M_0) = \qP(\overline{M_m})
$$
i.e. blowup first the subspace
$$
\qP(M_m/M_{m-1}) \subset \qP(\overline{M_m}) \subset \qP(M_m)
$$
then the strict transform of $\qP(M_m/M_{m-2})$ etc. There
is a natural isomorphism
$$
U(\sigma) \rightarrow q^{-1}(U).
$$
\end{Proposition}

\begin{Remark}\label{sean5.9.1}
When $\sigma = [M]$, \eqref{sean5.9} is immediate from \eqref{sean5.8}. As this
is the only case of \eqref{sean5.9} that we will need, we omit the
proof, which in any case is analogous to (and simpler than)
that of \eqref{sean5.5} and \eqref{sean5.5.2}.
\eqref{sean5.9} can also be deduced from
the claim on \cite[pg 168]{Faltings01} that for any
$[N] \in Y$ the natural map $\bS_Y \rightarrow \qP(N)$
is a composition of blowups with smooth centers (which
Faltings describes).
\end{Remark}

Now fix $\cF$ as in \eqref{sean1.3}.

\begin{Lemma}\label{sean5.10}\label{finitestable}
Let $Z \subset N$ be a subset
with $|Z| = r+1$. There is a unique stable lattice
$[\Lambda_Z] \in [\cF]$
such that the limits $f_Z^{\Lambda}$ are generic (i.e.
any $r$ of them is an $R$-basis).
In particular, there are finitely many stable lattices and $\Stab$ is finite.
If we reorder so that $Z = \{0,1,\dots,r\}$ and express
$$
f_0 = z^{a_1}p_1 f_1 + \dots z^{a_r} p_r f_r
$$
with $a_i \in \bZ$, $p_i \in R^*$, then
$\Lambda_Z = R z^{a_1} f_1 + \dots R z^{a_r} f_r$.

\end{Lemma}

\begin{proof} It's clear
that for $\Lambda_Z$ as given, the limits $\cF^{\Lambda_Z}$ are
in general position, so $\Lambda_Z$ is stable.

For uniqueness,
assume the limits $\cF^{\Lambda}_Z$ are in general position.
Then $f_i^{\Lambda}$, $r \geq i \geq 1$
are an $R$-basis of $\Lambda$. Define $b_i \in \bZ$ by
$z^{b_i} f_i = f_i^{\Lambda}$. Scaling $\Lambda$ we may
assume $b_i \geq a_i$, with equality for some $r \geq i \geq 1$.
Thus $f_0 = f_0^{\Lambda}$. Then
$b_i = a_i$, for all $i$, for otherwise $f_0^{\oLambda}$ will
be in the span of some proper subset of the $f_i^{\oLambda}$,
$r \geq i \geq 1$. So $\Lambda = \Lambda_Z$.
\end{proof}

\begin{Notation}\label{sean5.11}
For a subset $I \subset N$
let $V_I \subset V_K$ be the vector subspace spanned
by $f_i, i \in I$, and let $V^I := V/V_I$. For
each lattice $M \subset V$, let $M^I$ be its image in $V^I$,
i.e.
$M^I := M/M \cap V_I$.
\end{Notation}

Let $Y \subset [\cF]$ be a finite convex collection, containing
$\Stab$.
One checks immediately that the collection of equivalence
classes
$$
Y^I :=\{[M^I]\}_{[M] \in Y}
$$
is convex.

\begin{ETheorem}\label{sean5.12}\textsc{Definition--Lemma.}
Let $\bB_i \subset \bS_Y$
be the subfunctor of compatible quotients such that
$q_M$ vanishes on $f_i^M$ for all $[M] \in Y$.

Then $\bB_i$ is the Deligne scheme for $Y^{\{i\}}$.
$\bB_i \subset \bS_Y$ is non-singular and is the closure of
the hyperplane on the generic fibre
$$
\{f_i =0\} \subset \qP(V_K) \subset \bS_Y.
$$
\end{ETheorem}

\begin{proof} Clearly $M \cap V_I = f_i^M R$,
so its clear $\bB_i = \bS_{Y^{\{i\}}}$. The rest now
follows from \eqref{sean1.14}.
\end{proof}

\begin{Proposition}\label{sean5.13}
Let $[M] \in [\cF]$ be
a maximal lattice for a $k$ point of
$$
\cap_{i \in I} \bB_i \subset \bS_Y.
$$
Then the limits $f_i^{M}, i \in I$ are independent over $R$,
i.e. they generate an $R$ direct summand of $M$ of rank $|I|$.
\end{Proposition}

\begin{proof} We consider the corresponding simplex of
maximal lattices
$$
zM = M_0 \subset M_1 \dots \subset M_m = M.
$$
For each $m \geq s \geq 1$, let $I_s \subset I$ be those
$i$ so that $f_i^{M} \in M_s \setminus M_{s-1}$. Clearly
it is enough to show that the images of
$f_i^M$, $i \in I_s$, in $M_s/M_{s-1}$ are linearly independent.
By scaling (which allows us to move any of the $M_i$ to the $M_1$
position) it is enough to consider $s=1$, and show that
the images of $f_{i}$, $i \in I_1$ are linearly independent
in $M/zM$. Suppose not. Choose a minimal set whose images
are linearly dependent, which after reordering we may assume
are $f_0,f_1,\dots,f_p$. Further reordering we may assume
$$
f_1^{M}, f_2^M, \dots, f_r^M
$$
are an $R$-basis of $M$, or equivalently, their images given
a basis of $M/zM$. Renaming the $f_i$ we can assume
$f_i^M = f_i$. Now consider the unique expression
\begin{equation}\label{sean5.13.1}
f_0 = \sum_{i =1}^{r} z^{a_i} p_i f_i \cooltag
\end{equation}
with $a_i \in \bZ$, $p_i \in R^*$. By construction $a_i \geq 0$.
Since the images of $f_0,\dots,f_p$ in $M/zM$ are a minimal
linearly dependent set, it follows that $a_i \geq 1$ for $i > p$,
and $a_i =0$ for $p \geq i \geq 1$. Now let
$$
\Lambda := R f_1 + \dots R f_p + R z^{a_{p+1}} f_{p+1} +\dots
R z^{a_r} f_{r}.
$$

Note $f_i^{M_1} = f_i$ for $p \geq i \geq 0$, so by
assumption $q_{M_1}$ vanishes on these~$f_i$.
$z^{a_t} f_t \in zM_k = M_0 \subset M_1$ for $t \geq p+1$. Thus
by (5.2), $q_{M_1}$ vanishes on these $z^{a_t} f_t$. Thus
$q_{M_1}$ vanishes on $\Lambda$. But by equation \eqref{sean5.13.1}
and \eqref{sean5.10}, $\Lambda$ is stable. In particular $[\Lambda] \in Y$.
Clearly $\Lambda \subset M_1$, but $\Lambda \not \subset M_0$.
So the vanishing of $q_{M_1}|_{\Lambda}$ violates \eqref{sean5.2}.
\end{proof}

\begin{Theorem}\label{sean5.14}
Let $Y \subset [\cF]$ be a finite
convex set containing all the stable lattices.

For any subset $I \subset N$, the scheme theoretic intersection
$$
\bigcap_{i \in I} \bB_i \subset \bS_Y
$$
is non-singular and (empty or) codimension $|I|$
and represents
the Deligne functor $\bS_{Y^I}$. The divisor $S_Y + \bB \subset \bS_Y$
has normal crossings.
\end{Theorem}

\begin{proof} By \eqref{sean1.14} it's enough to show the intersection
represents
the Deligne functor.

Let $\cI$ be the intersection.
It is obvious from the definitions and \eqref{sean5.12} that
$\cI$ represents the subfunctor of $\bS_Y$ of
compatible quotients $q_M$ which vanish on $f_i^M$, $i \in M$,
while $\bS_{Y^I}$ is the subfunctor where $q_M$ vanishes on
$M \cap V_I$. Since $f_i^M \in M \cap V_I$, clearly
$$
\bS_{Y^I} \subset \cI
$$
is a closed subscheme. Since by \eqref{sean1.14}, $\bS_{Y^I}$ is
flat over $R$, it's enough to show that
they have the same special fibres. And so it is enough
to show the subfunctors agree on $T = \Spec(B)$ for
$B$ a local $k = R/zR$ algebra, with residue field $k$.
By maximal for such a $T$-point, we mean maximal for the
closed point. We take
a family of compatible quotients vanishing on all
$f_i^M$ and show they actually vanish on $V_I \cap M$.
By \eqref{sean5.2} it's enough to consider maximal $M$. But by \eqref{sean5.13}
$f_i^{M}$, $i \in I$ are independent over $R$, so in fact
they give an $R$-basis of $V_I \cap M$.
\end{proof}

\section{Global Generation}\label{sean6}
Here we complete the proof of \eqref{sean1.4}.

\begin{Review}\textsc{Bubbling.}
We choose an increasing sequence of finite convex subsets
$Y_i \subset \cB$ whose union is the membrane $[\cF]$ --
the existence for example follows from \eqref{sean4.3.2}.
We assume $\Stab \subset Y_i$.
We have the natural forgetful maps
$$
p_{i,j}: \bS_{Y_i} \rightarrow \bS_{Y_j},\quad i \geq j.
$$
By \eqref{sean5.14}, we have the vector bundle $\Omega^1_p(\log \bB)$ on $\bS_{Y_i}$,
see \ref{sean9}.
\end{Review}

\begin{Proposition}\label{sean6.2}
Given a closed point $x \in \bS_{Y_j}$
for all $i$ sufficiently large there is a closed point
$y \in \bS_{Y_i} \setminus \bB$
so that $p_{i,j}(y) = x${\rm:}
\end{Proposition}

\bigskip
\epsfxsize=0.8\textwidth
\centerline{\epsfbox{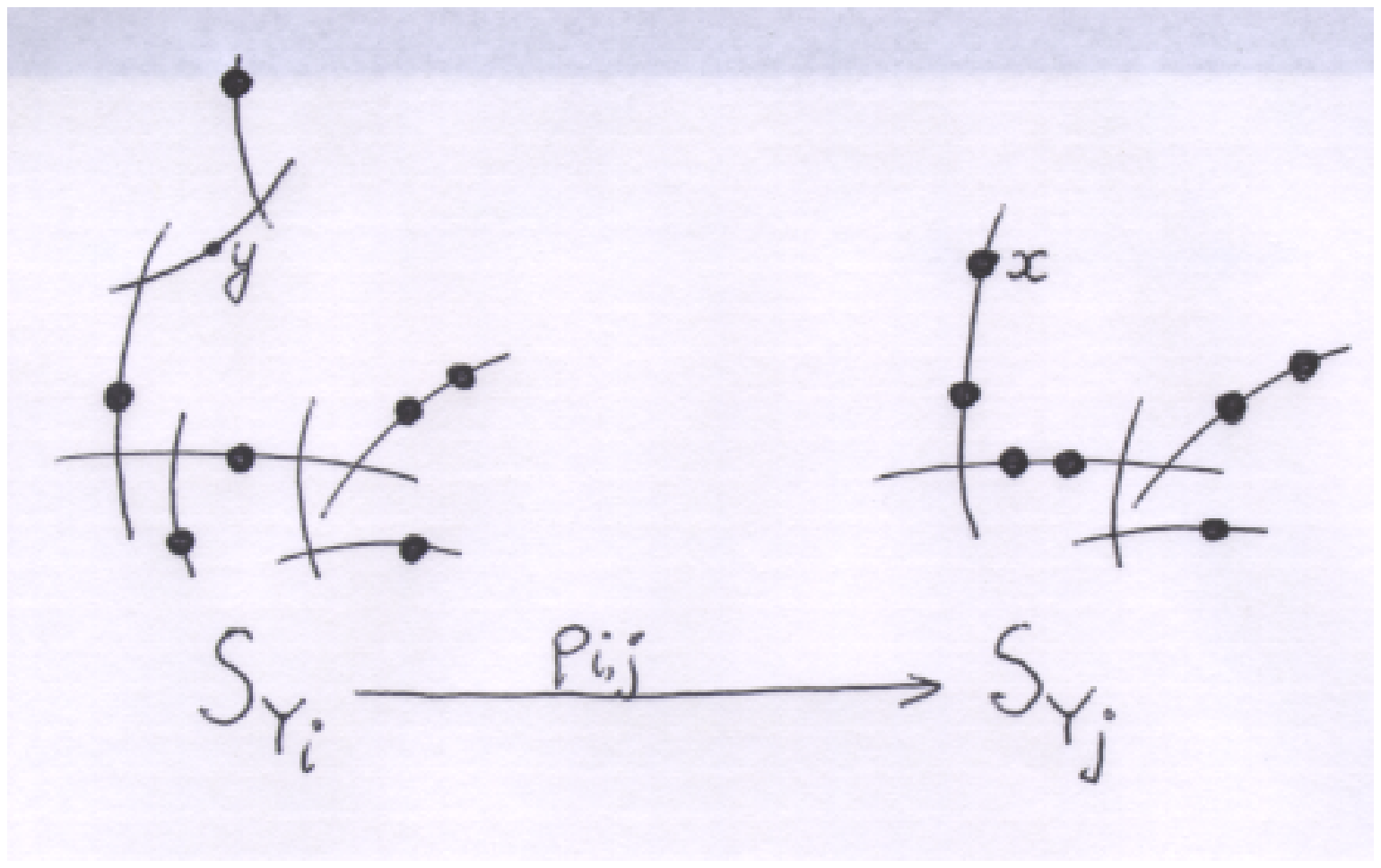}}
\bigskip

\begin{proof} Say $x$ lies on $\tqP_{Y_j}(\oM)$, $[M] \in Y_j$. Choose
$i$ sufficiently big so that
$$
[M + z^{-1}f^M R] \in Y_i
$$
for all $f \in \cF$. Clearly $p_{i,j}(\tqP_{Y_i}(\oM)) = \tqP_{Y_j}(\oM)$.
It remains to show that $\tqP_{Y_j}(\oM)$ is disjoint from
$\bB\subset \bS_{Y_j}$. Suppose that $q\in \tqP_{Y_j}(\oM)\cap\bB_k$.
Let $N=M + z^{-1}f_k^M R$.
Clearly $z^{-1}f_k^{M} = f_k^N$, so by definition of $\bB_k$,
$q_N$ vanishes on $z^{-1}f_k^{M}$. But $q_N$ vanishes
on $M$ by definition of a  maximal lattice. So $q_N =0$, a
contradiction.
\end{proof}

\begin{Review}\textsc{Barycentric Subdivision Trick.}
Next we introduce a convenient operation:
Let $R' = k[[z^{1/m}]]$ and $\Spec(R') \rightarrow \Spec(R)$
the associated finite map. Let $M_v$ be a collection
of lattices in $V_K$, and $Y$ their convex hull. Let
$M_v' := M_v \otimes_R R'$, and let $Y'$ be their
convex hull.
\end{Review}

\begin{Proposition}\label{sean6.3}
There is a commutative
diagram
$$
\begin{CD}
\bS_{Y'} @>{b}>> \bS_Y \\
@V{p'}VV            @V{p}VV \\
\Spec(R') @>>> \Spec(R)
\end{CD}
$$
with all arrows proper. If $m\ge r$ then given a $k$-point $y\in \bS_Y$
there
is a $k$-point $y'\in\bS_{Y'}$ in its inverse image that
lies on a unique irreducible component of the special
fibre $S_{Y'}${\rm:}
\end{Proposition}

\bigskip
\epsfxsize=0.8\textwidth
\centerline{\epsfbox{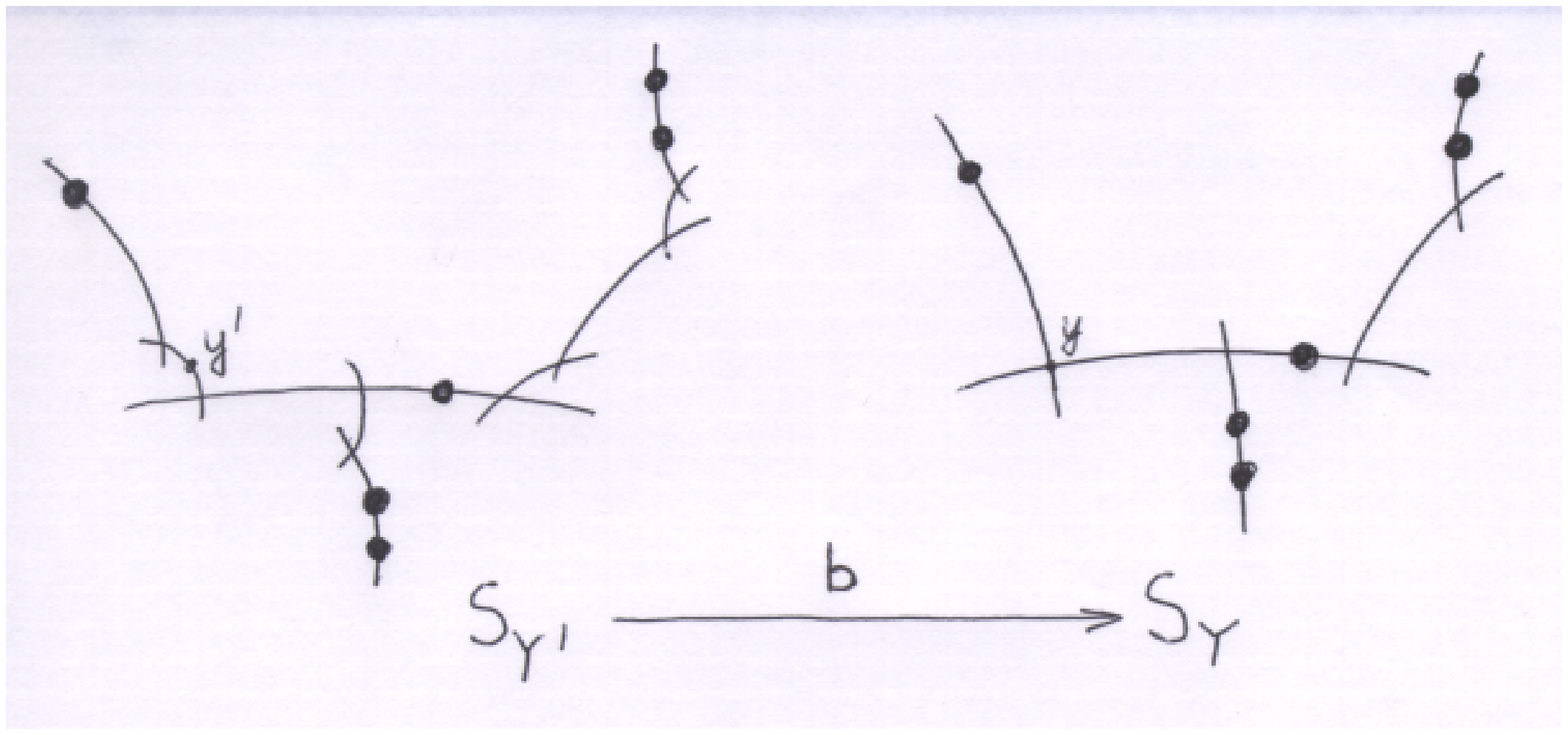}}
\bigskip

\begin{proof}
For any $R$-object $X$, we denote by $X'$ the base change to $R'$.
It is clear that $\bS_Y'$ represents
the functor $\bS_{\tilde Y}$ defined as in \eqref{sean1.14} but for the non-convex
collection $\tilde Y=\{M'\,|\,M\in Y\}$. Since $\tilde Y\subset Y'$,
there is a forgetful map $\bS_{Y'}\to\bS_Y'$ sending compatible
collections of quotients to compatible collections.
This implies the commutative diagram above.

Now let $m\ge r$, $\omega=z^{1/m}$, and let
$y\in \bS_Y$ be a $k$-point with the simplex of maximal lattices
$$\sigma=\{M_0=zM_k\subset M_1\subset\ldots\subset M_k\}.$$
By \eqref{sean5.2.3}, $y$ is determined by a collection of hyperplanes
$H_i \subset M_i/M_{i-1}$ which do not contain any
$\Res_{\sigma} M$ for $M \in Y$.
Let $N$ be the $R'$ lattice
$$N=M_1'+\omega M_2'+\omega^2 M_3'+\ldots+\omega^{k-1} M_k'.$$
Observe: $\omega^i M_i' \subset \omega N$, $ k \geq i \geq 0$.
($M_0' = \omega^{m-k} \omega^k M_k'$ and for $k > 0$ the
inclusion is clear). Thus we have a map
\begin{equation}\label{sean6.3.1}
\begin{CD}
\mathop\bigoplus\limits_{i=1}^{k} M_i'/M_{i-1}'
@>{\oplus \cdot \omega^{i-1}}>> N/\omega N.
\end{CD}\cooltag
\end{equation}
\eqref{sean6.3.1} is clearly surjective, thus it is an isomorphism,
as domain and range are $r$-dimensional $k$-vector spaces.
By the injectivity of the map we have
\begin{equation}\label{sean6.3.2}
\omega^{i-1} M_i' \cap \omega N = \omega^{i-1} M_{i-1}. \cooltag
\end{equation}
Now let $y' \in \qP(N/zN)$ be given by any hyperplane
$H' \subset N/zN$ which restricts to $H_i$ on $M_i/M_{i-1}$
under \eqref{sean6.3.1}. Its enough to show $y' \in U([N])$ for
then clearly $y'$
is send to $y$ and $N$ is the only  maximal lattice of $y'$ (or equivalently, by
\eqref{sean5.4}, $y'$ lies on a unique irreducible component of the special
fibre $S_{Y'}$).
 By \eqref{sean5.2.3} it's enough to show $H'$
does not contain $\Res_{[N]}[\Lambda]$ for
$[\Lambda] \in Y'$, and by \eqref{sean4.3.1} its enough to check
this for $[\Lambda] = [M']$, $[M] \in Y$.

We can assume
$M = M^{M_k}$, $M \subset M_i$, $M \not \subset M_{i-1}$. Then
$\Res_{\sigma}[M] = (M + M_i)/M_{i-1}$ and it
follows from \eqref{sean6.3.2} that $(M')^N = \omega^{i-1} M'$ thus we have
$$
(\cdot \omega^{i-1})(\Res_{\sigma} [M]) = \Res_{[N]} [M'].
$$
by \eqref{sean6.3.1}.
Thus $H'$ does not contain $\Res_{[N]} [M']$, since $H_i\not\supset\Res_{\sigma} [M]$.
\end{proof}

\begin{Theorem}\label{sean6.4}
Under the natural maps $p_{i,j}:\,\bS_{Y_i} \rightarrow \bS_{Y_j}$ for $j \leq i$,
there is an induced isomorphism of vector bundles
$$
p_{i,j}^*(\Omega^1_p(\log \bB)) \rightarrow \Omega^1_p(\log \bB).
$$
Furthermore the global sections
$$\dlog(f/g),\quad f,g\in\cF$$
generate $\Omega^1_p(\log \bB)$ globally. In particular
$\omega_p(\bB)$ is globally generated.
\end{Theorem}

\begin{proof} Since the sections are pulled back from
$\bS_{Y_j}$, the last remark will imply the first. Furthermore,
to prove the given sections generate at some point $y \in \bS_{Y_j}$,
it is enough to prove they generate at some inverse image point
on~$\bS_{Y_i}$. Thus by \eqref{sean6.2} it's enough to prove they generate
at a point $y \in \bS_i \setminus \bB$.

By \eqref{sean6.3}, we have the proper (generically finite) map
$$
\bS_{Y'_i} \rightarrow \bS_{Y_i}.
$$
Clearly the $\dlog(f/g)$ pullback to the analogous
forms on the domain, so by \eqref{sean6.3} we may assume $y$ has a unique maximal
lattice $M$. Now by
\eqref{sean5.9.1} the natural
map
$$
\bS_{Y_i} \rightarrow \qP(M)
$$
is an isomorphism of $y \in U([M])$
onto an open
subset of $\qP(M)$ which misses all of the hyperplanes
$$
f_i^M =0.
$$
Note $\cF^M$ contains an $R$-basis of $M$. Reordering,
say $f_k^M$, $r \geq k \geq 1$ give such a basis. Then it's
enough
to show the (regular) forms
$$
\dlog(f_k^M/f_1^M), r \geq k \geq 2
$$
give trivialisation of the ordinary cotangent bundle
over the open set in question, which is obvious.
\end{proof}

\begin{Review}\textsc{Minimal Model.}
By \eqref{sean6.4} the line bundle $\omega_p(\bB)$ is globally
generated. We consider the $p$-relative minimal
model, $\pi: \bS_Y \rightarrow \obS$, i.e.
\begin{equation}\label{sean6.5}
\obS :=\Proj \mathop\bigoplus\limits_m p_*(\omega_p(\bB)^{\otimes m}).
\cooltag
\end{equation}
Note by \eqref{sean6.4} that $\obS$ is independent of $Y$. Let
$\pi_*(\bB_i) =: \bB_i \subset \obS$.
\end{Review}

\begin{Theorem}\label{sean6.8}
Let $\Spec(R) \rightarrow \oX(r,n)$
be the unique extension of the map which sends the generic
fibre to
$$
[\bP^{r-1},L_1 + \dots L_n] \in \bP^0(r,n)/\PGL_r =
X(r,n)\subset\oX(r,n).$$
The pullback of the universal visible contour family
$(\cS,\cB)$, \eqref{sean2.7} is $(\obS,\bB)$.
\end{Theorem}

\begin{proof} By \eqref{sean6.4} we have a natural surjection
$$
V_n \otimes \cO_{\bS_Y} \twoheadrightarrow \Omega^1_p(\log \bB)
$$
inducing a regular map $\bS_Y \rightarrow G(r-1,n-1)$, which
on the general fibre is Kapranov's visible contour embedding
of $\bP^{r-1}$ given by the bundle of log forms with
poles on the $n$ general hyperplanes. This induces
a regular map $\bS_Y \rightarrow \cS$ where
$\cS \rightarrow \Spec(R)$ is the pullback of the visible
contour family, \eqref{sean2.7}.  $\omega_p(\bB)$ is pulled back from a
relatively very ample line bundle (the Pl\"ucker polarisation)
on $\cS$, by \eqref{sean2.12}. Thus $\bS_Y \rightarrow \cS$ factors through
a finite map
$$
\obS \rightarrow \cS.
$$
The map is birational, an isomorphism on the generic fibre.
By \eqref{sean2.8}, $\cS$ is normal. Thus it is an
isomorphism.
\end{proof}

\section{Bubble Space}

\begin{Review}
Here we prove \eqref{sean1.17}.
In \ref{sean6} we choose finite convex $Y\subset [\cF]$ containing
$\Stab$. Though there is a canonical choice, namely the convex hull
of $\Stab$, more esthetic is to take the infinite set
$[\cF]$. Let $Y_i$ be an increasing sequence of finite convex
subsets, containing $\Stab$, with union $[\cF]$.
Call $Y_i$ {\it full} along the simplex $\sigma \subset Y_i$
if $\Star_{[M]} Y_i = \Star_{[M]} [\cF]$ for all
$[M] \in \sigma$. It's clear that if $Y_i$ is full along~$\sigma$, so is $Y_j$ for $j \geq i$.
Let $U_i \subset Y_i$ be the union of all $U_{Y_i}(\sigma)$
such that $Y_i$ is full along~$\sigma$.
The next lemma shows that we may view
$U_i$ as an increasing sequence of open sets.
We define $\bS = \cup_{i} U_i.$
\end{Review}

\begin{Lemma}\label{sean6.11}
If $Y_i$ is full along $\sigma$ then $p_{ji}^{-1}(U_{Y_i}(\sigma)) \subset U_{Y_j}(\sigma)$.
Moreover, $p_{ji}^{-1}(U_i) \subset U_j$, and the map
$
p_{ji}^{-1}(U_i) \rightarrow U_i
$
is an isomorphism.
\end{Lemma}

\begin{proof}
Assume that $Y_i$ is full along $\sigma$.
Take $x \in p_{ji}^{-1}(U_{Y_j}(\sigma))$ and a maximal lattice
$[N] \in Y_j$ for $x$. Maximal
lattices form a simplex, so $[N]$ is adjacent
to a lattice in $\sigma$ and therefore
$[N] \in Y_i$ because $Y_i$ is full along $\sigma$.
Now it's clear $[N]$ is maximal for $p_{ji}(x)$,
so $[N] \in \sigma$. Thus $x \in U_{Y_i}(\sigma)\subset U_i$.

To show that $p_{ji}^{-1}(U_i) \rightarrow U_i$
is an isomorphism it suffices to check that
$p_{ji}^{-1}(U_{Y_i}(\sigma)) \rightarrow U_{Y_i}(\sigma)$ is an isomorphism
for any $\sigma\subset Y_i$. This map is
is proper and birational, and domain and range are non-singular,
so to show it is an isomorphism,
it's enough to check there are no-exceptional divisors, and
so to check that each irreducible component of the special fibre
of the domain maps onto an irreducible component of the special
fibre for the image. These components are the (appropriate
open subsets of the)
$\qP_{Y_j}(\oM)$, $[M] \in \sigma$, and its obvious
that $\qP_{Y_j}(\oM) \twoheadrightarrow \qP_{Y_i}(\oM)$.
\end{proof}

\begin{Theorem}\label{sean6.13}
$\bS$ is non-singular, and locally
of finite type. Its
special fibre, $S_{\infty}$, has smooth projective
irreducible components and
normal crossings. Let $\bB_i \subset \bS$ be the
hyperplane $f_i = 0$ of the generic fibre.
$\bB_i \subset \bS$ is closed and disjoint from $S_{\infty}$.
$\bB = \sum \bB_i$ has normal crossings.
In particular there is a natural
vector bundle $\Omega^1_{\bS/R}(\log \bB)$ whose determinant
is $\omega_{\bS/R}(\bB)$. The bundle is globally generated.

There are natural surjective
maps $p_i: \bS \rightarrow \bS_{Y_i}$ for all $i$, and
natural isomorphisms
$$
\begin{CD}
p_i^{-1}(\Omega^1_p(\log \bB)) = \Omega^1_{\bS/R}(\log \bB) \\
H^0(\bS_{Y_i},\Omega^1_{\bS_i/R}(\log \bB))
\rightarrow H^0(\bS,\Omega^1_{\bS/R}(\log \bB)). \\
\end{CD}
$$

The differential forms $\dlog(f/g)$, $f,g \in \cF$ define a
natural inclusion
$$
V_n \subset H^0(S_{\infty},\Omega^1_{S_{\infty}}).
$$
The sections generate the bundle. The associate map
$S_{\infty} \rightarrow G(r-1,n-1)$ factors through
$S_{Y_i}$ for all $i$, and its image is the special
fibre of the pullback of the family $\cS \rightarrow \oX(r,n)$
for the map $\Spec(R) \rightarrow \oX//H$ given by~$\cF$.
\end{Theorem}

\begin{proof}
Arguing as in \eqref{sean6.2},
one can show that the special fibre of $U_i \rightarrow \Spec(R)$
is disjoint from $\bB \subset \bS_{Y_i}$.
It follows that $S_\infty$ is disjoint from $\bB$.

We check that $\bS \rightarrow \bS_{Y_i}$ is surjective for
all $i$. The rest then follows easily from \eqref{sean6.4} and \eqref{sean6.8}.
Take $[M] \in Y_i.$
It's clear from the definitions that
$\tqP(\oM) \subset \bS_{Y_j}$ surjects onto
$\tqP(\oM) \subset \bS_{Y_i}$ for $ j \geq i$. Moreover,
there are only finitely many simplicies of $[\cF]$ that
contain $[M]$, by \eqref{sean4.2.1}, and for $j$ large $Y_j$ will contain
them all, from which it follows that $\tqP(\oM) \subset U_j$.
Thus the image of $\bS \rightarrow Y_i$ will contain
$\tqP(\oM)$. The union of the $\tqP(\oM)$ is the full special
fibre so $\bS \rightarrow Y_i$ is surjective.
\end{proof}



\begin{Review}
The Deligne functor for $[\cF]$ is not represented by a
scheme. However,
$\bS$ represents a natural subfunctor. In particular, $\bS$ is independent
on the choice of a sequence $Y_i$:
\end{Review}

\begin{Theorem}\label{sean6.17}
$\bS$ represents the subfunctor
of the Deligne functor for $[\cF]$ a $T$-valued
point of which is a collection of compatible quotients such
that each $k'$-point of $T$ admits a maximal lattice.
\end{Theorem}

\begin{proof} Take a $k$-point of $U_i$. By \eqref{sean6.11} any
lattice $[M] \in Y_i$ maximal in $Y_i$ for $k$ will be
maximal in $[\cF]$. It follows that $U_i$ is a subfunctor
of the functor in the statement, and thus $\bS$ is a subfunctor.
For the other direction its enough to consider $T$ the
spectrum of a local ring, with residue field $k'$. Consider
a $T$-point of the subfunctor in the statement.
Note that
in the proof of \eqref{sean5.2.3} the only place finiteness of $Y$ is
used is to establish the existence of a maximal lattice
which here we assume.
So the $k'$ point has a simplex of maximal
lattices, $\sigma$, satisfying \eqref{sean5.2.3}. For $i$ large, $Y_i$
will be full along $\sigma$, and now it is clear that the
quotients define a $T$-point of $U_i$, and thus of $\bS$.
\end{proof}

\section{Limit Variety}\label{sean7}

\begin{Review}
The matroid decomposition corresponding to the fiber of the visible
contour family can be readily
obtained from the power series, as we now describe. From the
matroid decomposition one can describe the fibre
using \cite[5.3]{Lafforgue03}.
We assume the reader is familiar with the general theory of variation
of GIT quotient, VGIT. See e.g. \cite{DH}. We note
in particular that $\Delta(r,n)$ parameterizes $\PGL_r$-ample
line bundles on $\bP(r,n)$ with non-empty semi-stable
locus.
\end{Review}

\begin{Definition}\label{sean6.14}
Call a polarisation
$L \in \Delta(r,n)$ on $\bP(r,n)$ {\it generic} if
there are no strictly semi-stable points.
\end{Definition}

\begin{Review}
For $[M] \in \cB$, let $C_M$ be the configuration of
limiting hyperplanes
$$\{f^{\oM} = 0\} \subset \qP(\oM),\quad f\in\cF$$
and let $P_M\subset\Delta(r,n)$ be the matroid polytope of $C_M$ (see~\eqref{sean3.1}).
\end{Review}

\begin{Lemma}\label{sean6.15}
If $L \in \Delta(r,n)$ is a generic
polarisation then there is a unique $[M] \in [\cF]$
so that the configuration $C_M$ is $L$-stable. $[M]$ is
GIT stable.
\end{Lemma}

\begin{proof} Let $Q$ be the GIT quotient of $\bP(r,n)$
given by $L$. $Q$ is a fine moduli space for $L$-stable
configurations and carries a universal family, a smooth
\'etale locally trivial $\bP^{r-1}$ bundle.
$\cF$ gives a $K$-point of $Q$, which extends uniquely to
an $R$-point. The pullback of the universal family will
be trivial over $R$ (as $R$ is Henselian), and so
$\qP(M)$ for some lattice $M \subset V_K$. It's clear the limit
configuration (given by the image in $Q$ of the closed
point of $R$) is equal to $C_M$. As $C_M$ has no automorphisms
it follows that this configuration contains $r$ hyperplanes
in general position, and so $[M] \in [\cF]$. The proof shows
that if $C_N,C_M$ are both $L$-stable, then the rational map
$\qP(N) \rightarrow \qP(M)$ is a regular isomorphism (either
is the pullback of the universal family over $Q$), which
implies $[N] = [M]$. $C_M$ has trivial automorphism
group, so the final remark holds by definition.
\end{proof}

\begin{Theorem}\label{sean6.16}
Let $x \in\oX(r,n)$  be
the limit point for the one parameter family given by $\cF$.
Assume $x$ belongs to a stratum given by the matroid decomposition $\uP$.
Then the maximal dimensional
polytopes of $\uP$ are precisely the $P_{M}$ for which
$C_{M}$ has no automorphisms, i.e.~$M$ is GIT-stable.
\end{Theorem}

\begin{proof}
By \cite{Kapranov93} and the theory of VGIT,
the matroid decomposition is
obtained as follows: GIT equivalence determines a polyhedral
decomposition of $\Delta(r,n)$. Chambers (interiors of maximal
dimensional polytopes in the decomposition) correspond to
polarisations with no strictly semi-stable points. For each
such chamber, there is a corresponding GIT quotient, which is
a fine moduli space for configurations of hyperplanes
stable for this polarisation. The one parameter family has
a unique limit in each such quotient, and in particular
associated to each chamber we have a limiting configuration.
Associated to
the configuration is its matroid polytope, and the polytopes
obtained in this way are precisely the facets of $\uP$.
Now by \eqref{sean6.14}, if $L$ is a polarisation in a chamber, then
the limiting configuration is $C_{M_L}$ for a unique
$[M_L]\in [\cF]$. Conversely, if we take $[M] \in [\cF]$
so that $C_{M}$ has no automorphisms, the polytope
$P_{M} \subset \Delta(r,n)$ is maximal dimensional,
see \cite[1.11]{Lafforgue03}. General $L \in P_{M}$
will be generic (in the sense of \eqref{sean6.14}) and it's clear
that $C_M$ is $L$-stable, so $M = M_L$.
\end{proof}

\begin{Review}\textsc{Stratification.}
The membrane $[\cF]$ is by \eqref{sean4.4.1} a union of apartments. We
have an apartment for each $I \subset N$, $|I| = r$, and thus
for each vertex of $\Delta(r,n)$. We stratify $[\cF]$ by apartments --
with one stratum for each collection of verticies -- those
points which lie in these, but no other, apartments. It follows
easily from \eqref{sean6.16} that the stratum is non-empty iff the
collection are the verticies of some $P \in \uP$, in which case
the stratum consists of those $[M] \in [\cF]$ (or
more generally, rational points of the realization, see \eqref{sean4.5})
with $P_M = P$. The dimension of the stratum is the codimension
of the polytope $P$ in $\Delta(r,n)$. We write $[\cF]^k$ for
the union of $k$-strata. Note $[\cF]^0 $ is precisely the
union of GIT stable lattices.

It is easy to describe the stratification in terms
of the power series $\cF$. Its enough to describe it in one apartment,
say $[f_1,\dots,f_r]$:

For any $a_1,\ldots,a_r$, let $S(a_1,\ldots,a_r)$ be the stratification of
$[f_1,\dots,f_r]$ by cones spanned by rays $R_0,\dots R_r$
where
$$
R_i= [z^{a_1}f_1,\ldots,z^{a_{i-1}}f_{i-1},
z^{a_{i}+p}f_{i}, z^{a_{i+1}}f_{i+1},\ldots,z^{a_r}f_r] \quad  p\ge0.
$$
\end{Review}

\begin{Lemma}\label{sean7.5}
Let $f_i=p_i^1z^{a^1_i}f_1+\ldots+p_i^rz^{a^r_i}f_r$ for $i=r+1,\ldots,n$,
where $p_i^k \in R^*$.
Then the stratification of $[f_1,\dots,f_r]$ is defined by
intersections of the $S(a^1_i,\ldots,a^r_i)$,
see the picture for $r =3$.
\end{Lemma}

\begin{proof} Immediate from the definitions \end{proof}

\bigskip
\epsfxsize=0.8\textwidth
\centerline{\epsfbox{pict5.ps}}
\bigskip

Now consider the limit pair $(S,B)$. The irreducible components
correspond to $[\cF]^{0}$, for $[M] \in [\cF]^{0}$ the corresponding
components is the log canonical model of $\qP(\oM) \setminus C_{[M]}$ --
the complement to the union of limiting hyperplanes.
We note by \cite[5.3]{Lafforgue03} that a collection of components
have a common point of intersection iff the corresponding matroids
have a common face, which is iff the points in $[\cF]^0$ all
lie on the boundary of the corresponding stratum. In particular,
if they have non-empty intersection, they all lie in a single
apartment.

\begin{Review}\textsc{Lines.} From now on  we assume that $r=3$.
\end{Review}

\begin{Lemma}\label{sean7.1}
A configuration of lines in
$\bP^2$ has trivial automorphism group
iff it contains $4$ lines in linear general
position. A configuration has non-trivial automorphism group
iff there is a point in the configuration which is in the
complement of at most one of the lines. In this case the
automorphism group is positive dimensional.
\end{Lemma}

\begin{proof} This is easy linear algebra.
\end{proof}

\begin{Remark}
By \eqref{sean7.1}, stable is the same as GIT stable if $r=3$.
This fails in higher dimensions: the configuration
of planes
$$x_1=0,\ x_2=0,\ x_3=0,\ x_4=0,\quad x_1+x_2+x_3=0,\quad x_2+x_3+x_4=0$$
in $\bP^3$ is GIT stable but not stable.
\end{Remark}

\begin{Lemma}\label{sean7.2}
Let $C$ be a stable configuration of
lines indexed by $N$. Let $\tS \rightarrow \bP^2$ be the blowup
of all points of multiplicity at least $3$. Let $B \subset \tS$ be
the reduced inverse image of the lines. then $K_{\tS} + B$ is
ample and
$$
\bP^2 \setminus C \subset \tS
$$
is the log canonical compactification except in one case:
If there are two points $a,b \in L$ on a line $L$ of $C$ such that
any other line of $C$ meets $L$ in either $a$ or $b$. In this case
(the strict transform of) $L \subset \tS$ is a $(-1)$-curve, and the
blowdown is $\bP^1 \times \bP^1$,
$$
\bP^2 \setminus C \subset \bP^1 \times \bP^1
$$
is the log canonical model, with boundary a union of fibers for
the two rulings.
\end{Lemma}

We refer to this exceptional case as a {\it special
stable configuration}:

\bigskip
\epsfxsize=0.4\textwidth
\centerline{\epsfbox{pict6.ps}}
\bigskip

\begin{proof} We induct on $n$. When $n =4$ then $\tS = \bP^2$
and the result is obvious. So we assume $n > 4$. If the configuration
is special, the result is clear, so we assume it is not.
By \eqref{sean7.1} we can drop a line, $M$,
so the resulting configuration $C'$ is stable. If $C'$ is special,
with special line $L$, then since $C$ is not special, it follow
that if we instead drop $L$, the resulting configuration is
stable, and non-special. So we may assume $C'$ is not special.
Add primes to the notation to indicate analogous objects for $C'$.
We have $q: \tS \rightarrow \tS'$, the blowup along the points of $M$
where $C$ has multiplicity exactly $3$ (note $\tS' \rightarrow \bP^2$
is an isomorphism around these points). Thus
$$
K_{\tS} + B = q^*(K_{\tS'} + B') + M
$$
(where we use the same symbol for a curve and for
its strict transform).
It's clear $K_{\tS} + B$ is $q$-ample. As $K_{\tS'} + B'$ is
ample, the only curve  on which $K_{\tS} + B$ can have non-positive
intersection is $M$. But $(K_{\tS} + B) \cdot M > 0$ by
adjunction, since $C$ is not special. It follows that $K_{\tS} + B$
is ample.
\end{proof}

\begin{Review}\label{sean7.4}\textsc{The Limit Surface.}
Now we describe the limit pair
$(S,B)$ precisely. The irreducible components are smooth, and
described by \eqref{sean7.2}. We write $S_{M}$ for the component
corresponding to $[M] \in\Stab$.
Unbounded $1$-strata -- rays in some apartment,
correspond to irreducible components of $B$, bounded $1$-strata
correpond to irreducible components of $\Sing(S)$. For each
$[M] \in\Stab$, the $1$-strata which bound $[M]$ correspond to
boundary components of $S_{M}$ (components of the complement
to $\qP(\oM) \setminus C_{M}$). $S_{M}$ and $S_{N}$ have
one dimensional intersection iff $[M],[N]$ is the boundary of
a $1$-stratum, and in that case they are glued along the
corresponding boundary component (a copy of $\bP^1$). Triple
points of $S$ (points on three or more components) correspond
to bounded $2$-strata. The local analytic singularities of
$(S,B)$ are described by the following:
\end{Review}

\begin{Theorem}\label{sean7.6}
Let $p \in S$ be a point where the
pair $(S,B)$ fails to have normal crossings. There are two possibilities
for the germ
of $(S,B)$ in an analytic neighborhood  $p \in U_p$:
\begin{equation}\label{sean7.6.1}
U_p =\langle e_1,e_2\rangle\cup\langle e_2,e_3\rangle\cup \langle e_{3},e_{4}\rangle\subset \bA^{4}, \cooltag
\end{equation}
$B \cap U_p$ is the union of $\langle e_1\rangle$ and $\langle e_4\rangle$, and these are components
of a single $B_i$.
\begin{equation}\label{sean7.6.2}
U_p=\langle e_1,e_2 \rangle\cup \langle e_2,e_3 \rangle\cup\ldots\cup \langle e_n,e_1 \rangle\subset \bA^n,
\cooltag
\end{equation}
$n=3,4,5,6$. $B \cap U_p = \emptyset$.

Here $e_1,\dots,e_n$ are coordinate axes in $\bA^n$, and
$\langle\cdot\rangle$ is the linear span.
\end{Theorem}

\begin{proof}
It is simple to classify bounded $2$-strata using \eqref{sean7.5}. Then the
glueing among components is described by \cite[5.3]{Lafforgue03}.
\end{proof}

\section{The bundle of relative log differentials}\label{sean9}
Here we recall a general construction which we will use
at various points throughout the paper:

Let
$$
p:(\cS,\cB) \rightarrow C,\quad \cB = \sum_{i=1}^{n} \cB_i,
$$
be a pair
of a non-singular variety with normal crossing divisor,
semi-stable over the curve $C$, in a neighborhood of
$0 \in C$, i.e. $(\cS,F + \cB)$ has normal crossings, where
$F$ is the fibre over $0$. We assume the general fibre is
projective, but not necessarily the special fibre.
We define the bundle of
relative log differentials $\Omega^1_{p}(\log \cB)$ by
the exact sequence
\begin{equation}\label{sean9.1}
\begin{CD}
0 @>>> \Omega^1_{C/k}(\log 0) \rightarrow
\Omega^1_{\cS/k}(\log F + \cB)
\rightarrow \Omega^1_p(\log \cB) \rightarrow 0
\end{CD}\cooltag
\end{equation}

Assume on the generic fibre $S$ that the restrictions
of the boundary components, $B_i,B_j$ are linearly
equivalent. Then we can choose a rational function
$f$ on $\cS$ so that
$$
(f) = \cB_i - \cB_j + E
$$
for $E$ supported on $F$.
Then $\dlog(f)$ gives a global section of
$\Omega^1_{\cS/k}(\log F + \cB)$. Note $f$ is unique up
to multiplication by a unit on $C \setminus 0$, and thus
the image of $\dlog(f)$ in $\Omega^1_p(\log \cB)$, which
we denote by $\dlog(\cB_i/\cB_j)$ is independent of $f$.

>From now on we assume that for the general fibre
$H^0(S,\Omega^1_S) =0$.  $\dlog(\cB_i/\cB_j)$
is now characterized as the unique section whose restriction
to the general fibre has residue $1$ along $B_i$, $-1$
along $B_j$ and is regular off of $B_i + B_j$.

In this way we obtain a canonical map
\begin{equation}\label{sean9.2}
V_n \rightarrow H^0(\cS,\Omega^1_p(\log \cB)) \cooltag
\end{equation}
($V_n$ the standard $k$-representation of the symmetric group
$S_n$) which is easily seen to be injective, e.g. by the
description of the residues on the general fibre.

The restriction
\begin{equation}\label{sean9.3}
\Omega^1_{S/k}(\log B) := \Omega^1_p(\log \cB)|_S \cooltag
\end{equation}
for $(S,B)$ the special fibre of $(\cS,\cB)$,
is canonically associated to $(S,B)$, i.e. is independent
of the smoothing. See e.g. \cite[\S3]{Friedman} or
\cite{KN} -- these authors treat
normal crossing
varieties without boundary, but the theory extends to
normal crossing pairs in an obvious way. Finally
there is a canonical residue map (e.g induced via \eqref{sean9.1} by
ordinary residues on $\cS$)
\begin{equation}\label{sean9.4}
\res : \Omega^1_{S_Y/k}(\log B)|_{B_i} \rightarrow \cO_{B_i} \cooltag
\end{equation}
together with \eqref{sean9.2} this gives a canonically split inclusion
\begin{equation}\label{sean9.5}
V_n \subset H^0(S,\Omega^1_S(\log B)) \cooltag.
\end{equation}

\begin{Definition}\label{sean9.6} Let $(S,B)$ be a normal
variety with boundary. Assume for an open subset
$i:U \subset S$ with complement of codimension at least two
that $U$ is non-singular and $B|_U$ has normal crossings.
Define
$$
\Omega^1_{S/k}(\log B) := i_*(\Omega^1_{U/k}(\log B|_U).
$$
\end{Definition}

\end{document}